\documentclass[12pt]{amsart}
\usepackage{amssymb}
\usepackage{graphicx}
\newcommand{\C}{\mathcal{C}}

\newtheorem{Def}{Definition}
\newtheorem{Lem}{Lemma}

\newtheorem{Thm}{Theorem}
\newtheorem{Cor}{Corollary}
\newtheorem{Rem}{Remark}
\newenvironment{Pf}{ Proof.}{\(\square\)}

\newtheorem{Exc}{Exercise}

\title[On compatible linear connections...]{On compatible linear connections of two-dimensional generalized Berwald manifolds}
\author{Cs. Vincze, T. Khoshdani, S. Mehdi Zadeh and M. Ol\'{a}h}
\address{Inst. of Math., Univ. of Debrecen \\
H-4010 Debrecen, P.O.Box 12 \\
Hungary}
\email{csvincze@science.unideb.hu}
\email{khoshdani@yahoo.com, zahira.m.2012@gmail.com, olma4000@gmail.com}
\begin{document}
\begin{abstract}In the paper we present results about generalized Berwald surfaces involving the intrinsic characterization, some topological obstructions for the base manifold and examples. 
\end{abstract}
\maketitle
{\centering\footnotesize {\emph{In memoriam to V. Wagner on the 75th anniversary of publishing his pioneering work about generalized Berwald manifolds.}} \par}
\footnotetext[1]{Cs. Vincze is supported by the EFOP-3.6.1-16-2016-00022 project. The project is co-financed by the European Union and the European Social Fund.}
\footnotetext[2]{T. Khoshdani and S. Mehdi Zadeh are supported by the Department of Mathematics of University of Mohaghegh Ardabili, Ardabil, Iran. M. Ol\'{a}h is supported by the University of Debrecen (Summer Grant 2018).}
\section*{Introduction}

The concept of generalized Berwald manifolds goes back to V. Wagner \cite{Wag1}. They are Finsler manifolds admitting linear connections such that the parallel transports preserve the Finslerian length of tangent vectors (compatibility condition). To express the compatible linear connection in terms of the canonical data of the Finsler manifold is the problem of the intrinsic characterization we are going to solve in case of two-dimensional generalized Berwald manifolds. The result is formulated in terms of linear inhomogeneous differential equations for the main scalar along the indicatrix curve (Subsection 2.1). As an application we prove that if a Landsberg surface is a generalized Berwald manifold then it must be a Berwald manifold (Subsection 2.2). Especially, we reproduce Wagner's original result in terms of the conventional setting of Finsler surfaces (Subsection 2.3) in honor of the 75th anniversary of publishing his pioneering work about generalized Berwald manifolds.

The technic of averaging is an alternative way to solve the problem of the characterization of compatible linear connections. By the fundamental result of the theory \cite{V5} such a linear connection must be metrical with respect to the averaged Riemannian metric given by integration of the Riemann-Finsler metric on the indicatrix hypersurfaces. Therefore the linear connection is uniquely determined by its torsion tensor. The torsion tensor has a special decomposition in 2D because of 
\begin{equation} 
\label{tor}
T(X,Y)=\left(X^1Y^2-X^2Y^1\right)\left(T_{12}^1 \frac{\partial}{\partial u^1}+T_{12}^2 \frac{\partial}{\partial u^2}\right)=\rho(X)Y-\rho(Y)X,
\end{equation}
where $\rho_1=T_{12}^2$ and $\rho_2=-T_{12}^1=T_{21}^1$. In higher dimensional spaces such a linear connection is called semi-symmetric. Using some previous results \cite{V6}, \cite{V9}, \cite{V10} and \cite{V11}, the torsion tensor of a semi-symmetric compatible linear connection can be expressed in terms of metrics and differential forms given by averaging independently of the dimension of the space. Especially the compatible linear connection must be of zero curvature in 2D unless the manifold is Riemannian, see \cite{V12}. Therefore we can conclude some topological obstructions as well due to the divergence representation of the Gauss curvature (Subsection 3.1). We prove, for example, that any compact generalized Berwald surface without boundary must have zero Euler characteristic. Therefore the Euclidean sphere does not carry such a geometric structure. Using the theory of closed Wagner manifolds, this means that the local conformal flatness of the Riemannian surfaces fails for (non-Riemannian) Finslerian ones (Subsection 3.2). We present some examples of non-Riemannian two-dimesnional generalized Berwald manifolds as well (Subsection 3.3).

\section{Notations and terminology}

Let $M$ be a differentiable manifold with local coordinates $u^1, \ldots, u^n.$ The induced coordinate system of the tangent manifold $TM$ consists of the functions $x^1, \ldots, x^n$ and $y^1, \ldots, y^n$. For any $v\in T_pM$, $x^i(v)=u^i(p)$ and $y^i(v)=v(u^i)$, where $i=1, \ldots, n$ and $\pi\colon TM\to M$ is the canonical projection.
\subsection{Finsler metrics} A Finsler metric is a continuous function $F\colon TM\to \mathbb{R}$ satisfying the following conditions:
\begin{itemize}
\item[(F1)] $\displaystyle{F}$ is smooth on the complement of the zero section (regularity),
\item[(F2)] $\displaystyle{F(tv)=tF(v)}$ for all $\displaystyle{t> 0}$ (positive homogenity),
\item[(F3)] the Hessian $\displaystyle{g_{ij}=\frac{\partial^2 E}{\partial y^i \partial y^j}}$, where $E=\frac{1}{2}F^2$ 
is positive definite at all nonzero elements $\displaystyle{v\in T_pM}$ (strong convexity).
\end{itemize}

The so-called \emph{Riemann-Finsler metric} $g$ is constituted by the components $g_{ij}$. It is defined on the complement of the zero section. The Riemann-Finsler metric makes each tangent space (except at the origin) a Riemannian manifold with standard canonical objects such as the {\emph {volume form}} $\displaystyle{d\mu=\sqrt{\det g_{ij}}\ dy^1\wedge \ldots \wedge dy^n}$,
the \emph {Liouville vector field} $\displaystyle{C:=y^1\partial /\partial y^1 +\ldots +y^n\partial / \partial y^n}$ together with its \emph{normalized dual form} $\displaystyle{l_i=\partial F/\partial y^i}$ with respect to the Riemann-Finsler metric and the {\emph {induced volume form}}
$$\mu=\sqrt{\det g_{ij}}\ \sum_{i=1}^n (-1)^{i-1} \frac{y^i}{F} dy^1\wedge\ldots\wedge dy^{i-1}\wedge dy^{i+1}\ldots \wedge dy^n$$
 on the indicatrix hypersurface $\displaystyle{\partial K_p:=F^{-1}(1)\cap T_pM\ \  (p\in M)}$. In what follows we summarize some basic notations. As a general reference of Finsler geometry see \cite{BSC} and \cite{M1}: $\displaystyle{g^{ij}=(g_{ij})^{-1}}$ denotes the inverse of the coefficient matrix of the Riemann-Finsler metric, the (lowered) first Cartan tensor is given by $\displaystyle{C_{ijk}=\frac{1}{2}\partial g_{ij}/\partial y^k }$  and $\displaystyle{\mathcal{C}^l_{ij}=g^{lk}\mathcal{C}_{ijk}}$. The first Cartan tensor is totally symmetric and $\displaystyle{y^k\C_{ijk}=0}$. Its semibasic trace is given by the quantities $\displaystyle{\mathcal{C}_i=g^{jk}\mathcal{C}_{ijk}}$ ($i, j, k=1, \ldots, n$).
Differentiating $\det g_{ij}$ as a composite function we have that
$$
\frac{\partial \det g_{rs}}{\partial y^i}=\frac{\partial D}{\ \partial m_{jk}} (M) \frac{\partial g_{jk}}{\partial y^i}=(-1)^{j+k}\det \left(M \ \textrm{ without its $j^{th}$ row and $k^{th}$ column}\right)\frac{\partial g_{jk}}{\partial y^i}=
$$
$$=(\det g_{rs}) g^{jk}\frac{\partial g_{jk}}{\partial y^i}, \ \ \textrm{where}\ \ M:=g_{ij}.$$
Therefore
\begin{equation}
\label{mainscalar1}
\frac{\partial \ln \sqrt{\det g_{rs}}}{\partial y^i}=\frac{1}{2} g^{jk}\frac{\partial g_{jk}}{\partial y^i}=g^{jk}\mathcal{C}_{ijk}=\mathcal{C}_i.
\end{equation}
The geodesic spray coefficients and the horizontal sections are
$$
G^l=\frac{1}{2}g^{lm}\left(y^k\frac{\partial^2 E}{\partial y^m\partial x^k}-\frac{\partial E}{\partial x^m}\right)
\ \ \textrm{and}\ \ X_i^h=\frac{\partial }{\partial x^i}-G^l_i\frac{\partial }{\partial y^l},\ \ \textrm{where}\ \ G_i^l=\frac{\partial G^l}{\partial y^i}.$$
The second Cartan tensor (Landsberg tensor) and the mixed curvature are given by
$$P^l_{ij}=\frac{1}{2} g^{lm}\left( X_i^h \left(g_{jm}\right)-G^k_{ij}g_{km}-G^k_{im} g_{jk}\right),\ \ \textrm{where}\ \ G_{ij}^l=\frac{\partial G_{i}^l}{\partial y^j}$$
and $\displaystyle{P_{ijk}^l=-G^{l}_{ijk}}$, where $\displaystyle{G_{ijk}^l=\frac{\partial G_{ij}^l}{\partial y^k}}$.
\begin{Lem}
\begin{equation}
\label{eq:4}
P^l_{ij}=-\frac{F}{2}l_m g^{kl}P_{ijk}^m
\end{equation}
\end{Lem}
\begin{Pf} Since
$$Fl_m=\frac{\partial E}{\partial y^m}, \ \ \frac{\partial E}{\partial y^m}G^m=\frac{1}{2}y^k\frac{\partial E}{\partial x^k}, \ \ g_{mi}G^m=\frac{1}{2}\left(y^k\frac{\partial^2 E}{\partial y^i\partial x^k}-\frac{\partial E}{\partial x^i}\right)$$
and
$$\frac{\partial} {\partial y^i}\left(\frac{\partial E}{\partial y^m}G^m\right)-g_{mi}G^m=\frac{\partial E}{\partial x^i}$$
we have
$$-Fl_mP^m_{ijk}=\frac{\partial E}{\partial y^m}G^m_{ijk}=$$
$$\frac{\partial}{\partial y^k}\left(\frac{\partial E}{\partial y^m}G^m_{ij}\right)-g_{mk}G_{ij}^m=\frac{\partial}{\partial y^k}\left(\frac{\partial}{\partial y^j} \left(\frac{\partial E}{\partial y^m}G^m_{i}\right)-g_{mj}G_i^m\right)-g_{mk}G_{ij}^m=$$
$$\frac{\partial}{\partial y^k}\left(\frac{\partial}{\partial y^j} \left(\frac{\partial} {\partial y^i}\left(\frac{\partial E}{\partial y^m}G^m\right)-g_{mi}G^m\right)-g_{mj}G^m_i\right)-g_{mk}G_{ij}^m=$$
$$\frac{\partial}{\partial y^k}\left(\frac{\partial}{\partial y^j} \left(\frac{\partial E}{\partial x^i}\right)-g_{mj}G^m_i\right)-g_{mk}G_{ij}^m=\frac{\partial}{\partial x^i}g_{jk}-2\C_{jmk}G_i^m-g_{mj}G_{ik}^m-g_{mk}G^m_{ij}=$$
$$2P_{ijk}=2g_{kl}P^l_{ij}\ \ \Rightarrow\ \ P^l_{ij}=-\frac{F}{2}l_m g^{kl}P_{ijk}^m$$
as was to be proved.
\end{Pf}

\subsection{Generalized Berwald manifolds} 

\begin{Def} A linear connection $\nabla$ on the base manifold $M$ is called \emph{compatible} to the Finslerian metric if the parallel transports with respect to $\nabla$ preserve the Finslerian length of tangent vectors. Finsler manifolds admitting compatible linear connections are called generalized Berwald manifolds.
\end{Def}

\begin{Cor}
\label{cor:1}
A linear connection  $\nabla$ on the base manifold $M$ is \emph{compatible} to the Finslerian metric function if and only if the induced horizontal distribution is conservative, i.e. the derivatives of the fundamental function $F$ vanish along the horizontal directions with respect to $\nabla$.
\end{Cor}

\begin{Pf} Suppose that the parallel transports with respect to $\nabla$ (a linear connection on the base manifold) preserve the Finslerian length of tangent vectors and let $X_t$ be a parallel vector field along the curve $c\colon [0,1]\to M$:
\begin{equation}
\label{paralleldiff}
(x^k\circ X_t)'={c^k}'\ \ \textrm{and}\ \ (y^k \circ X_t)'={X_t^k}'=-{c^i}'  X_t^j  \Gamma_{ij}^k\circ c
\end{equation}
because of the differential equation for parallel vector fields. If $F$ is the Finslerian fundamental function then
\begin{equation}
\label{eq:5}
(F \circ X_t)'=(x^k\circ X_t)'{\frac{\partial F}{\partial x^k}}\circ X_t+(y^k \circ X_t)'{\frac{\partial F}{\partial y^k}}\circ X_t
\end{equation}
and, by formula (\ref{paralleldiff}),
\begin{equation}
\label{eq:55}
(F \circ X_t)'={c^i}'\bigg(\frac{\partial F}{\partial x^i}-y^j {\Gamma}_{ij}^{k}\circ \pi \frac{\partial F}{\partial y^k}\bigg)\circ X_t.
\end{equation}
This means that the parallel transports with respect to $\nabla$ preserve the Finslerian length of tangent vectors (compatibility condition) if and only if
\begin{equation}
\label{cond1}
\frac{\partial F}{\partial x^i}-y^j {\Gamma}^k_{ij}\circ \pi \frac{\partial F}{\partial y^k}=0\ \  (i=1, \ldots,n),
\end{equation}
where the vector fields of type
\begin{equation}
\label{eq:6}
\frac{\partial}{\partial x^i}-y^j {\Gamma}^k_{ij}\circ \pi \frac{\partial}{\partial y^k}
\end{equation}
span the associated horizontal distribution belonging to $\nabla$.
\end{Pf}

\begin{Thm}
\label{heritage} \emph{\cite{V5}} If a linear connection on the base manifold is compatible with the Finslerian metric function then it must be metrical with respect to the averaged Riemannian metric
\begin{equation}
\label{averagemetric1}
\gamma_p (v,w):=\int_{\partial K_p} g(v, w)\, \mu=v^i w^j \int_{\partial K_p} g_{ij}\, \mu \ \ (v, w\in T_p M, p\in U).
\end{equation}
\end{Thm}

\subsection{Finsler surfaces} In case of Finsler surfaces it is typical to introduce the vector field
$$V:=\frac{\partial F}{\partial y^1}\frac{\partial}{\partial y^2}-\frac{\partial F}{\partial y^2}\frac{\partial}{\partial y^1}.$$
It is tangential to the indicatrix curve because of $VF=0$. Since three vertical vector fields must be linearly dependent in 2D,
$$$$
$$
0=\det \left(\begin{array}{ccc}
g\left(\frac{\partial}{\partial y^1}, \frac{\partial}{\partial y^1}\right)& g\left(\frac{\partial}{\partial y^1}, \frac{\partial}{\partial y^2}\right)&g\left(\frac{\partial}{\partial y^1}, C\right)\\
&&\\
g\left(\frac{\partial}{\partial y^2}, \frac{\partial}{\partial y^1}\right)& g\left(\frac{\partial}{\partial y^2}, \frac{\partial}{\partial y^2}\right)&g\left(\frac{\partial}{\partial y^2}, C\right)\\
&&\\
g\left(C, \frac{\partial}{\partial y^1}\right)& g\left(C, \frac{\partial}{\partial y^2}\right)&g\left(C,C \right)\\
&&\\
\end{array}\right)=\det \left(\begin{array}{ccc}
g_{11}& g_{12}& \partial E/\partial y^1\\
&&\\
g_{12}& g_{22}&\partial E/\partial y^2\\
&&\\
\partial E/\partial y^1& \partial E/\partial y^2& 2E\\
&&\\
\end{array}\right)=
$$
$$F^2 \det g_{ij}+2g_{12}\frac{\partial E}{\partial y^1} \frac{\partial E}{\partial y^2}-\left(\frac{\partial E}{\partial y^1}\right)^2 g_{22}-\left(\frac{\partial E}{\partial y^2}\right)^2g_{11}=F^2 \left(\det g_{ij}-g(V,V)\right).$$
$$$$
This means that $\displaystyle{0 \neq \det g_{ij}=g(V,V)}$ and, consequently, 
$$V_0:=\frac{1}{\sqrt{g(V,V)}}V, \ \ C_0:=\frac{1}{F}C, \ \ V_0^h:=V_0^i X_i^h=V_0^i \left(\frac{\partial }{\partial x^i}-G^l_i\frac{\partial }{\partial y^l}\right),\ \ S_0:=\frac{1}{F}S=\frac{y^i}{F} X_i^h$$ 
form  a local frame on the complement of the zero section in $\pi^{-1}(U)$. Such a collection of vector fields is called a \emph{Berwald frame} on the tangent manifold.
\begin{Def}
The main scalar of a Finsler surface is defined as $\displaystyle{\lambda:=V_0^jV_0^kV_0^l \mathcal{C}_{jkl}},$
where $\displaystyle{V_0=V/\sqrt{g(V,V)}}$ is the unit tangential vector field to the indicatrix curve.
\end{Def}
The vanishing of the main scalar implies that the surface is Riemannian and vice versa. The zero homogeneous version $I:=F\lambda$ of the main scalar is also frequently used in the literature \cite{Berwald1}, \cite{Berwald2}, \cite{H2} and \cite{M1}. Consider the vector field $\displaystyle{\mathcal{C}^k_{ij} \partial/\partial y^k}$. Since it is also tangential to the indicatrix surface it follows that
$$\mathcal{C}^k_{ij} \frac{\partial}{\partial y^k}=\mathcal{C}^l_{ij}g\left(V_0, \frac{\partial}{\partial y^l}\right)V_0,$$
where $\displaystyle{V_0=V/\sqrt{g(V,V)}}$ is the unit tangential vector field to the indicatrix curve. Therefore
$$\mathcal{C}^k_{ij}=\mathcal{C}^l_{ij} g_{lm}V_0^m V_0^k=V_0^m \mathcal{C}_{ijm} V_0^k \ \ \Rightarrow \ \ \mathcal{C}_{ijr}=V_0^m \mathcal{C}_{ijm} V_0^k g_{kr}.$$
Contracting by $g^{rj}$
\begin{equation}
\label{trace}
\mathcal{C}_i=V_0^j V_0^m \mathcal{C}_{ijm}.
\end{equation}
By formulas (\ref{mainscalar1}) and (\ref{trace}) we have that
\begin{equation}
\label{mainscalar2}
\lambda:=V_0^jV_0^kV_0^l \mathcal{C}_{jkl}=V_0^j \mathcal{C}_j=V_0 \left(\ln \sqrt{\det g_{rs}}.\right)
\end{equation}
In what follows we summarize some of the general formulas to express the surviving components of the Landsberg tensor, the mixed curvature tensor and the pairwise Lie-brackets of a Berwald frame (Cartan's permutation formulas) \cite{VV1}:
\begin{equation}
\label{wag01}
y^i V_0^j V_0^k P_{ijk}=y^iV_0^j V_0^k G^l_{ijk} g\left( V_0, \frac{\partial}{\partial y^l} \right)=0,
\end{equation}
$$V_0^i  V_0^j V_0^k P_{ijk}=-S(\lambda), \ \ V_0^i V_0^j V_0^k G^l_{ijk} g\left( V_0, \frac{\partial}{\partial y^l} \right)=V_0^h(\lambda)+V_0(S\lambda)$$
because of the homogenity properties; see \cite[Corollary 1.8.]{VV1} and \cite[Formula (24a)]{VV1}. E. Cartan's permutation formulas are 
\begin{equation}
\label{wag02}
[V_0, V_0^h]=-\frac{1}{F}S_0-\lambda V_0^h -S(\lambda) V_0, \ \ [S_0, V_0]=-\frac{1}{F}V_0^h, \ \ [V_0^h,S_0]=- \kappa V_0,
\end{equation}
where $\kappa$ is the only surviving coefficient of the curvature of the horizontal distribution \cite[Theorem 1.10]{VV1}. 
Let the indicatrix curve in $T_pM$ be parameterized as the integral curve of $V_0$:
$$V_0\circ c_p(\theta)=c_p'(\theta)\ \ \Rightarrow \ \ \lambda\circ c_p (\theta)=\left(\ln \sqrt{\det g_{rs}}\circ c_p\right)'(\theta).$$
It is called the \emph{central affine arcwise parametrization of the indicatrix curve}. The parameter $\theta$ is "the central affine length of the arc of the indicatrix" and the main scalar can be interpreted as its "central affine curvature"; for the citations see \cite{Wag1}.

\section{Two-dimensional generalized Berwald manifolds}

Let $\nabla$ be a linear connection on the base manifold $M$ and suppose that the parallel transports preserve the Finslerian length of tangent vectors (compatibility condition). By Corollary \ref{cor:1}, 
$$\frac{\partial E}{\partial x^i}-y^ m\Gamma_{im}^l\circ \pi  \frac{\partial E}{\partial y^l}=0 \ \ (i=1, 2).$$

\subsection{The comparison of $\nabla$ with the canonical horizontal distribution of the Finsler manifold} Using the canonical horizontal sections we can write that
$$y^m\Gamma_{im}^l\circ \pi \frac{\partial E}{\partial y^l}-G^{l}_i \frac{\partial E}{\partial y^l}=0.$$
Since the vertical vector fields are the linear combinations of $V$ and $C$, it follows that 
$$y^m\Gamma_{im}^l\circ \pi\frac{\partial }{\partial y^l}-G^{l}_i \frac{\partial }{\partial y^l}=f_i V+g_i C \ \ (i=1, 2);$$
the coefficients $f_1$, $f_2$ are positively homogeneous of degree one, $g_1$ and $g_2$ are positively homogeneous of degree zero. Since $VE=0$ but $CE=2E$, we have that $g_1=g_2=0$ and, consequently
\begin{equation}
\label{zero}
y^m\Gamma_{im}^l\circ \pi \frac{\partial }{\partial y^l}-G^{l}_i \frac{\partial }{\partial y^l}=f_i V \ \ \Rightarrow \ \ y^m\Gamma_{im}^k\circ \pi \frac{\partial }{\partial y^k}= G^{k}_i \frac{\partial }{\partial y^k}+f_i V \ \ (i=1, 2).
\end{equation}
To provide the linearity of the right hand side we should take the Lie brackets with the vertical coordinate vector fields  two times:
$$0=\left [ \left [y^m\Gamma_{im}^l\circ \pi \frac{\partial }{\partial y^l}, \frac{\partial }{\partial y^j} \right ], \frac{\partial }{\partial y^k}\right ]=\left [ \left [G^{l}_i \frac{\partial }{\partial y^l}, \frac{\partial }{\partial y^j}\right ], \frac{\partial }{\partial y^k}\right ]+\left [ \left [f_i V, \frac{\partial }{\partial y^j}\right ], \frac{\partial }{\partial y^k} \right ]=$$
$$G_{ijk}^l\frac{\partial }{\partial y^l}+f_i\left [ \left [V, \frac{\partial}{\partial y^j}\right ], \frac{\partial }{\partial y^k} \right ]-\frac{\partial f_i}{\partial y^j}\left [V, \frac{\partial}{\partial y^k}\right ]-\frac{\partial f_i}{\partial y^k}\left [V, \frac{\partial}{\partial y^j}\right ]+\frac{\partial^2 f_i}{\partial y^j \partial y^k} V=:W_{ijk},$$
where
$$\left[V, \frac{\partial}{\partial y^j}\right ]=\frac{\partial^2 F}{\partial y^j \partial y^2}\frac{\partial}{\partial y^1}-\frac{\partial^2 F}{\partial y^j \partial y^1}\frac{\partial}{\partial y^2},\ \ \left [ \left [V, \frac{\partial}{\partial y^j}\right ], \frac{\partial }{\partial y^k} \right ]=-\frac{\partial^3 F}{\partial y^j \partial y^k \partial y^2}\frac{\partial}{\partial y^1}+\frac{\partial^3 F}{\partial y^j \partial y^k \partial y^1}\frac{\partial}{\partial y^2}.$$
Since $\displaystyle{y^jW_{ijk}=y^kW_{ijk}=0}$ it is enough to investigate the quantity $W_i=V^jV^kW_{ijk}$. By some direct computations
$$V^j\frac{\partial^2 F}{\partial y^j \partial y^2}=V\left(\frac{\partial F}{\partial y^2}\right)=\frac{1}{F}V\left(F\frac{\partial F}{\partial y^2}\right)= \frac{1}{F}g\left(V, \frac{\partial}{\partial y^2}\right)$$
because of $VF=0$. On the other hand
$$V^jV^k\frac{\partial^3 F}{\partial y^j \partial y^k \partial y^2}=\frac{1}{F}V^k V \left(F \frac{\partial^2 F}{\partial y^k \partial y^2}\right)=\frac{1}{F}V^k V\left(g_{k2}-\frac{\partial F}{\partial y^k}\frac{\partial F}{\partial y^2}\right)=$$
$$\frac{1}{F}\left(2V^jV^k\mathcal{C}_{jk2}-V^kV\left(\frac{\partial F}{\partial y^k}\right)\frac{\partial F}{\partial y^2} \right)=\frac{1}{F}\left(2V^jV^k\mathcal{C}_{jk2}-\frac{1}{F}V^kV\left(F\frac{\partial F}{\partial y^k}\right)\frac{\partial F}{\partial y^2} \right)=$$
$$\frac{1}{F}\left(2V^jV^k\mathcal{C}_{jk2}-\frac{1}{F}g(V,V)\frac{\partial F}{\partial y^2} \right)$$
and, consequently,
$$W_i=V^jV^k G_{ijk}^l\frac{\partial }{\partial y^l}-\frac{2V(f_i)}{F}\left(g\left(V, \frac{\partial}{\partial y^2}\right) \frac{\partial}{\partial y^1}-g\left(V, \frac{\partial}{\partial y^1}\right)\frac{\partial}{\partial y^2}\right)-$$ 
\begin{equation}
\label{wdef}
\frac{f_i}{F}\left(\left(2V^jV^k \mathcal{C}_{jk2}-\frac{1}{F}g(V,V)\frac{\partial F}{\partial y^2} \right)\frac{\partial}{\partial y^1}-\left(2V^jV^k \mathcal{C}_{jk1}-\frac{1}{F}g(V,V)\frac{\partial F}{\partial y^1} \right)\frac{\partial}{\partial y^2}\right)+
\end{equation}
$$V^jV^k\frac{\partial^2 f_i}{\partial y^j \partial y^k} V.$$
The vanishing of $W_i$ is equivalent to
\begin{equation}
g(W_{i}, V_0)=0 \ \ \textrm{and}\ \ g(W_{i}, C_0)=0 \ \ (i=1, 2),
\end{equation}
where $V_0=V/\sqrt{g(V,V)}$ and $C_0=C/F$ are the normalized vector fields of the vertical Berwald frame. 

\subsubsection{The vanishing of the orthogonal term to the indicatrix} It follows that
$$0=g(W_{i}, C)=W_{i} E =FV^jV^kG_{ijk}^l\frac{\partial F}{\partial y^l}-2V(f_i)g(V,V)-2f_iV^jV^kV^l \mathcal{C}_{jkl}.$$
Therefore
\begin{equation}
\label{difeq01}
\frac{\alpha_i}{\sqrt{g(V,V)}}=\lambda f_i+ \left(V_0 f_i\right) \ \ (i=1,2),
\end{equation}
where $V_0=V/\sqrt{g(V,V)}$ is the unit tangential vector field to the indicatrix curve, $\lambda$ is the main scalar and 
$$\alpha_i=\frac{1}{2}FV_0^j V_0^k G_{ijk}^l\frac{\partial F}{\partial y^l}\stackrel{(\ref{eq:4})}{=} V_0^j V_0^k P_{ijk}.$$
Using that $\displaystyle{\det g_{ij}=g(V,V)}$, formula (\ref{mainscalar2}) says that
\begin{equation}
\label{difeq015}
\alpha_i=V_0\left( f_i \sqrt{g(V,V)}\right) \ \ (i=1,2).
\end{equation}
Let the indicatrix curve $c_p$ in $T_pM$ be parameterized as the integral curve of $V_0$. Evaluating along $c_p$ we have
\begin{equation}
\label{difeq02}
\alpha_i \circ c_p(\theta)=\left( f_i\circ c_p  \sqrt{g(V,V)}\circ c_p\right)'(\theta) \ \ (i=1,2)
\end{equation}
for any $p\in U$. Therefore
\begin{equation}
\label{difeq02hom}
\beta_i\circ c_p(t)=f_i\circ c_p(t) \sqrt{g(V,V)}\circ c_p(t)-f_i\circ c_p(0) \sqrt{g(V,V)}\circ c_p(0),
\end{equation}
where $\displaystyle{\beta_i \colon \pi^{-1}(U)\to \mathbb{R}}$ $(i=1, 2)$ are the $1$-homogeneous extensions of the functions defined by
\begin{equation}
\label{difeq035}
\beta_i \circ c_p(t)=\int_0^t \alpha_i\circ c_p(\theta)\, d\theta \ \ (i=1, 2)
\end{equation}
along the central affine arcwise parametrization of the indicatrix curve. We can write that 
\begin{equation}
\label{difeq03}
f_i \circ c_p(t)=\frac{1}{\sqrt{g(V,V)}\circ c_p(t)} \left( \beta_i\circ c_p(t)+k_i(p)\right) \ \ (i=1, 2)
\end{equation}
for some constants $k_i(p)$ ($i=1, 2$) depending only on the position. 

\subsubsection{The vanishing of the tangential term to the indicatrix} It follows that
$$0 =g(W_{i}, V)= V^j V^k G^l_{ijk} g\left( V , \frac{\partial}{\partial y^l}\right)
 - \frac{2 f_{i}}{F} \left(V^j V^k \mathcal{C}_{jk2} g\left( V, \frac{\partial}{\partial y^1}\right) -  V^j V^k \mathcal{C}_{jk1} g\left( V, \frac{\partial}{\partial y^2} \right)\right)+$$
$$ \frac{ f_{i}}{F^2}  g (V, V) \left(\frac{\partial F}{\partial y^2} g\left( V , \frac{\partial}{\partial y^1} \right) -  \frac{\partial F}{\partial y^1} g \left( V, \frac{\partial}{\partial y^2} \right)\right)
 + V^j V^k \frac{\partial ^2 f_{i}}{ \partial y^j \partial y^k} g(V,V)=V^j V^k G^l_{ijk} g \left( V, \frac{\partial}{\partial y^l} \right) $$
$$+ V^j V^k \frac{\partial ^2 f_{i}}{ \partial y^j \partial y^k} g(V,V)-\frac{ f_{i}}{F^2}g^2 (V, V)-\frac{2 f_{i}}{F} \left(V^j V^k \mathcal{C}_{jk2} g\left( V, \frac{\partial}{\partial y^1} \right) -  V^j V^k \mathcal{C}_{jk1} g\left( V, \frac{\partial}{\partial y^2}\right)\right),$$
where
$$V^j V^k \mathcal{C}_{jk2} g\left( V, \frac{\partial}{\partial y^1} \right) -  V^j V^k \mathcal{C}_{jk1} g\left( V, \frac{\partial}{\partial y^2}\right)=0$$
because the vector field 
$$Z:=g\left( V, \frac{\partial}{\partial y^1} \right)\frac{\partial}{\partial y^2} -  g\left( V, \frac{\partial}{\partial y^2}\right)\frac{\partial}{\partial y^1}$$
is parallel to $C$, i.e. $\displaystyle{g(V,Z)=0}$. Therefore 
$$0=V^j V^k G^l_{ijk} g\left( V, \frac{\partial}{\partial y^l} \right) + V^j V^k \frac{\partial ^2 f_{i}}{ \partial y^j \partial y^k} g(V,V)-\frac{ f_{i}}{F^2}g^2 (V, V)$$
and, consequently, 
\begin{equation}
0= V_0^j V_0^k G^l_{ijk} g\left( V_0, \frac{\partial}{\partial y^l} \right) + V_0^j V_0^k \frac{\partial ^2 f_{i}}{ \partial y^j \partial y^k} \sqrt{g(V,V)} - \frac{ f_i}{ F^2}  \sqrt{ g (V, V)}. \label{formula3}
\end{equation}
\begin{Lem} \label{lem:2} If $g$ is a positively homogeneous function of degree $k$, then 
\begin{equation}
V_0(V^k_0) \frac{\partial g}{ \partial y^k} = - \lambda V_0 (g) - k\frac{g}{ F^2}.\label{formula1}
\end{equation}
Especially,
\begin{equation}
V_0^j V_0^k \frac{\partial ^2 f_{i}}{ \partial y^j \partial y^k} = V_0 (V_0 f_i) +\lambda V_0 (f_i) +\frac{f_i}{ F^2}.
\end{equation}
\end{Lem}
\begin{Pf}
Let $c_p$ be the parametrization of the indicatrix curve in $T_pM$ as the integral curve of $V_0$, i.e. $V_0 \circ c_p = c_p'$. Differentiating equation
\begin{eqnarray}
1= g_{c_p} ( V_0 \circ c_p , V_0 \circ c_p ) = g_{ij} \circ c_p  (c_{p}^{i})'  (c_p^{j}) '\label{T1}
\end{eqnarray}
we have that $\displaystyle{0 = 2 g_{ij} \circ c_p (c^{i}_p)'' (c^{j}_p)' + 2 \mathcal{C}_{ijk} \circ c_p (c_p^{i})' (c_p^{j})' (c_p^{k})'}$ and, consequently, 
\begin{equation}
g_{c_p}(V_0\circ c_p, c_p'')=g_{c_p} ( c_p' , c_p '') = - \mathcal{C}_{ijk} \circ c_p (c_p^i)' (c_p^j)' (c_p^k)' = - \left(V_0^iV_0^jV_0^k\mathcal{C}_{ijk}\right) \circ c_p= - \lambda \circ c_p.\label{curve1}
\end{equation}
Differentiating equation 
\begin{eqnarray}
 0= g_{c_p} ( C \circ c_p , V_0 \circ c_p) = g_{ij} \circ c_p (c^i_p) (c_p^j)'\label{T2}
\end{eqnarray}
we have that 
\begin{equation}
0= 2 \mathcal{C}_{ijk} \circ c_p (c^i_p) (c_p^j)'(c_p^k)'+ g_{ij} \circ c_p (c_p^{i})' (c_p^{j})'+ g_{ij} \circ c_p (c_p)^i (c_p^j)''.\label{curve2}
\end{equation}
Taking into account that $\displaystyle{\mathcal{C}_{ijk} \circ c_p (c^i_p) (c_p^j)'(c_p^k)' =  \mathcal{C}_{ijk} \circ c_p  (y^i \circ c_p) (c_p^j)'(c_p^k)'= 0}$, 
$$g_{ij} \circ c_p (c_p^{i})' (c_p^{j})'=g_{c_p} (c_p' , c_p')=1 \ \ \textrm{and} \ \ g_{ij} \circ c_p (c_p)^i (c_p^j)''= g_{c_p} ( C\circ c_p , c^{''}_p), $$
it follows that
\begin{equation}
g_{c_p} ( C_0 \circ c_p , c^{''}_p ) =- \frac{1}{F\circ c_p} \label{curve3},
\end{equation}
where $C_0:=C/F$ is the normalized Liouville vector field. From (\ref{curve1}) and (\ref{curve3})
\begin{equation}
c^{''}_p  = - (\lambda V_0) \circ c_p - \frac{1}{F\circ c_p}  C_0 \circ c_p .\label{curve4}
\end{equation}
This means that 
$$
\left( V_0(V^k_0) \frac{\partial g}{ \partial y^k} \right) \circ c_p = (V^k_0 \circ c_p )'  \frac{\partial g}{ \partial y^k} \circ c_p  = (c_p^k)'' \frac{\partial g}{ \partial y^k} \circ c_p \overset{(\ref{curve4})}{=}$$
$$ -\left((\lambda V^k_0) \circ c_p + \frac{1}{F\circ c_p}  C^k_0 \circ c_p  \right)\frac{\partial g}{ \partial y^k} \circ c_p = - (\lambda V_0 g) \circ c_p - \frac{1}{F^2\circ c_p} (C g) \circ c_p, $$
where $\displaystyle{Cg=kg}$ because of the homogenity. Note that the terms $\displaystyle{V_0(V^k_0) \partial g/ \partial y^k}$, $\displaystyle{\lambda V^k_0 g}$ and $\displaystyle{g/F^2}$ are of the same degree of homogenity, i.e. they are homogeneous of degree $k-2$. Therefore the equality along the indicatrix curve implies (\ref{formula1}). Especially,
$$V_0 (V_0 f_i)= V_0^j V_0^k \frac{\partial ^2 f_{i}}{ \partial y^j \partial y^k}+V_0(V^k_0) \frac{\partial f_i}{ \partial y^k} =V_0^j V_0^k \frac{\partial ^2 f_{i}}{ \partial y^j \partial y^k} -\lambda V_0 (f_i) -\frac{f_i}{ F^2}$$
as was to be proved.
\end{Pf}

Using Lemma \ref{lem:2} we can write formula (\ref{formula3}) into the form 
\begin{equation}
\label{formula001}
0= \omega_i + \left(V_0(V_0f_i)\right) \sqrt{g(V,V)}+\lambda V_0(f_i) \sqrt{g(V,V)}, 
\end{equation}
where
$$\omega_i=V_0^j V_0^k G^l_{ijk} g\left( V_0, \frac{\partial}{\partial y^l} \right) \ \ (i=1, 2).$$
By formula (\ref{mainscalar2})
\begin{equation}
0= \omega_i + \left(V_0(V_0f_i)\right) \sqrt{g(V,V)}+V_0(f_i)V_0 \left(\sqrt{g(V,V)}\right) \ \ (i=1, 2)
\end{equation}
because of $\displaystyle{\det g_{ij}=g(V,V)}$. Therefore
\begin{equation}
0= \omega_i + V_0\left((V_0 f_i) \sqrt{g(V,V)} \right),  \label{formulamain01} 
\end{equation}
\begin{equation}
0= \omega_i + V_0\left(V_0\left (f_i \sqrt{g(V,V)}\right)-f_i V_0\left(\sqrt{g(V,V)}\right)\right) , \label{formulamain02} 
\end{equation}
\begin{equation}
0= \omega_i + V_0\left(V_0\left (f_i \sqrt{g(V,V)}\right)- \lambda f_i \sqrt{g(V,V)}\right), \label{formulamain03}
\end{equation}
\begin{equation}
0= \omega_i + V_0\left(V_0\left (f_i \sqrt{g(V,V)}\right)\right)- V_0(\lambda) f_i \sqrt{g(V,V)}-\lambda V_0\left( f_i \sqrt{g(V,V)}\right) \ \ (i=1, 2). \label{formulamain04}
\end{equation}
By formula (\ref{difeq015})
\begin{equation}
0= \omega_i + V_0\left(\alpha_i\right)- V_0(\lambda) f_i \sqrt{g(V,V)}-\lambda \alpha_i \ \ (i=1, 2).\label{formulamain045}
\end{equation}
Evaluating formula (\ref{formulamain045}) along $c_p$
\begin{equation}
\omega_i\circ c_p (t) + (\alpha_i\circ c_p)'(t)= \left(\beta_i\circ c_p(t)+k_i(p)\right)(\lambda\circ c_p)'(t)+\lambda\circ c_p(t) \alpha_i\circ c_p(t) \ \ (i=1, 2) \label{formulamain05}
\end{equation}
because of (\ref{difeq035}) and (\ref{difeq03}). The constants $k_1(p)$ and $k_2(p)$ of integration can be expressed by (\ref{formulamain05}) provided that $\displaystyle{\lambda \circ c_p}$ is not a constant function:
$$
k_i(p)=\frac{ \gamma_i\circ c_p(s)- \beta_i\circ c_p(s) \lambda \circ c_p(s)+\alpha_i\circ c_p(s)-\alpha_i\circ c_p(0)}{\lambda \circ c_p(s)-\lambda \circ c_p(0)},$$
where 
$$\gamma_i\circ c_p(s)=\int_0^s \omega_i\circ c_p(t)\, dt\ \ (i=1, 2)$$
and the parameter $s\in \mathbb{R}$ is choosen such that $\displaystyle{\lambda \circ c_p(s)-\lambda \circ c_p(0)\neq 0}$. Otherwise the function $\displaystyle{\lambda\circ c_p}$ is constant. Since $\det g_{ij}$ attains its extremals along the indicatrix curve, formula (\ref{mainscalar2}) shows that $\displaystyle{\lambda\circ c_p}$ is identically zero and the indicatrix is a quadratic curve in $T_pM$. The quadratic indicatrix curve of a (connected) generalized Berwald manifold at a single point implies that the indicatrices are quadratic curves at any point and we have a Riemannian surface. Indeed, the parallel transports induced by the compatible linear connection take a quadratic curve into quadratic curves\footnote{Non-Riemannian Finsler surfaces with main scalar depending only on the position must be singular; see Berwald's original list  \cite[Formulas 118 I-III]{Berwald1}, see also \cite{Berwald2} and \cite{VV}.}.

\begin{Thm} The compatible linear connection of a non-Riemannian connected generalized Berwald surface must be of the form
\begin{equation}
\Gamma_{ij}^1\circ \pi= G^{1}_{ij} -\frac{\partial f_i}{\partial y^j}\frac{\partial F}{\partial y^2}-f_i \frac{\partial^2 F}{\partial y^j \partial y^2}, \ \ \Gamma_{ij}^2\circ \pi= G^{2}_{ij} +\frac{\partial f_i}{\partial y^j}\frac{\partial F}{\partial y^1}+f_i \frac{\partial^2 F}{\partial y^j \partial y^1}\ \ (i, j=1, 2),
\end{equation}
where the functions $f_1$, $f_2$ are given by 
\begin{equation}
f_i \circ c_p(t)=\frac{1}{\sqrt{g(V,V)}\circ c_p(t)} \left( \int_0^t \alpha_i\circ c_p(\theta)\, d\theta+k_i(p)\right) \ \ (i=1, 2)
\end{equation}
and the integration constants satisfy equations
\begin{equation}
\omega_i\circ c_p (t) + (\alpha_i\circ c_p)'(t)= \left(\int_0^t \alpha_i\circ c_p(\theta)\, d\theta+k_i(p)\right)(\lambda\circ c_p)'(t)+\lambda\circ c_p(t) \alpha_i\circ c_p(t) \ \ (i=1, 2)
\end{equation}
for any $p\in M$.
\end{Thm}

\begin{Pf} Equations for the functions $f_1$ and $f_2$ imply that $g(W_i,C)=0$ because of subsection 2.1.1. Equations for the integration constants imply that $\displaystyle{g(W_i, V_0)=0}$ because of subsection 2.1.2. Therefore $W_i=0$ and we have a generalized Berwald surface. The explicite formulas for the coefficients of the linear connection preserving the Finslerian length of tangent vectors are 
\begin{equation}
\Gamma_{ij}^1\circ \pi= G^{1}_{ij} +\frac{\partial f_i}{\partial y^j}V^1-f_i \frac{\partial^2 F}{\partial y^j \partial y^2}, \ \ \Gamma_{ij}^2\circ \pi= G^{2}_{ij} +\frac{\partial f_i}{\partial y^j}V^2+f_i \frac{\partial^2 F}{\partial y^j \partial y^1}\ \ (i, j=1, 2),
\end{equation}
because of formula \eqref{zero}.
\end{Pf}

\begin{Cor} The compatible linear connection of a generalized Berwalds surface is uniquely determined. 
\end{Cor}

\subsection{An application: Landsberg and generalized Berwald surfaces}

\begin{Def}
A Finsler manifold is called a Landsberg manifold if the Landsberg tensor of the canonical horizontal distribution vanishes. The Berwald manifolds are defined by the vanishing of the mixed curvature tensor of the canonical horizontal distribution.
\end{Def}

Formula (\ref{eq:4}) implies that any Berwald manifold is a Landsberg manifold. The converse of this statement is the famous Unicorn problem in Finsler geometry \cite {B}.

\begin{Thm}
A connected non-Riemannian two-dimensional generalized Berwald surface is a Landsberg surface if and only if it is a Berwald surface. 
\end{Thm}

\begin{Pf}
Suppose that we have a connected non-Riemannian two-dimensional generalized Berwald manifold such that the Landsberg tensor vanishes, i.e. $\alpha_i=0$ ($i=1, 2$). Then (\ref{difeq03}) implies that
\begin{equation}
f_i \sqrt{g(V,V)} = k_i(p) F\label{difeq06}
\end{equation}
for any point $p\in M$. On the other hand 
\begin{equation}
\omega_i - V_0(\lambda) f_i \sqrt{g(V,V)} =0 \label{difeq04}
\end{equation}
due to (\ref{formulamain045}). Contracting by $y^i$ 
\begin{eqnarray}
V_0(\lambda) y^i k_i(p) =0.\label{difeq05}
\end{eqnarray}
If there exists a point $p \in M$ such that $\displaystyle{k_1^2 (p) + k_2^2(p) \neq 0.}$ then $\displaystyle{y^1 k_1(p) + y^2 k_2(p)=0}$ is an equation of a line in $T_p M$. Therefore, there are at most two positions along $\partial K_p$ such that $\displaystyle{v^1 k_1(p) + v^2 k_2(p) =0}$. Otherwise $V_0 (v) \lambda =0$ because of (\ref{difeq05}). 
A continuity argument says that $ V_0(v) \lambda =0 $ for any $v \in T_p M $, i.e. $\lambda$ is constant along $c_p$. Since $\det g_{ij}$ attains its extremals along the indicatrix curve, formula (\ref{mainscalar2}) shows that $\lambda \circ c_p =0$. This means that the indicatrix is a quadratic curve in $T_p M$. The quadratic indicatrix curve of a (connected) generalized Berwald manifold at a single point implies that the indicatrices are quadratic curves at any point due to the compatible linear connection and the induced linear mapping between the tangent spaces. This is a contradiction because the generalized Berwald surface is non-Riemannian. Otherwise $\displaystyle{k_1(p) = k_2(p) =0}$ for any $p \in M$, i.e. $f_i=0$ ($i=1, 2$) and the compatible linear connection must be the canonical one. Therefore we have a Berwald manifold. 
\end{Pf}

\subsection{Wagner's equations} To present Wagner's equations in \cite{Wag1} we need the following simple observation: 
$$H_i=0 \ \ (i=1, 2) \ \ \textrm{if and only if}\ \ y^iH_i=0 \ \ \textrm{and}\ \ V_0^i H_i=0$$
because of 
$$
\det \left(\begin{array}{cc}
y^1& y^2\\
V_0^1& V_0^2\\
\end{array}\right)=y^1V_0^2-y^2V_0^1=\frac{F}{\sqrt{g(V,V)}}\neq 0.
$$ 
Contracting (\ref{formulamain045}) by $y^i$
$$0=y^i V_0\left(\alpha_i\right)- V_0(\lambda) y^i f_i \sqrt{g(V,V)},$$
where $\displaystyle{y^iV_0 (\alpha_i)=V_0(y^i\alpha_i)-V_0^i\alpha_i\stackrel{(\ref{wag01})}{=}S(\lambda)}$ and, consequently,
\begin{equation}
\label{wag04}
S(\lambda)= V_0(\lambda) y^i f_i \sqrt{g(V,V)}.
\end{equation}
Contracting (\ref{formulamain045}) by $V_0^i$
$$0\stackrel{(\ref{wag01})}=V_0^h\lambda+V_0(S\lambda)+V_0^iV_0\left(\alpha_i\right)-V_0(\lambda) V_0^i f_i \sqrt{g(V,V)}+\lambda S(\lambda),$$
where 
$$V_0^iV_0\left(\alpha_i\right)=V_0\left(V_0^i \alpha_i\right)-V_0(V_0^i)\alpha_i \stackrel{(\ref{wag01})}=-V_0(S\lambda)-V_0(V_0^i)\alpha_i.$$
Since $V_0(V_0^i)\circ c_p=\left(c_p^i\right)''$ it follows, by formula (\ref{curve4}), that 
\begin{equation}
\label{wag015}
V_0(V_0^i)=-\lambda V_0^i-\frac{y^i}{F^2}
\end{equation}
due to the $-1$-homogeneous extension. Therefore
$$V_0^iV_0\left(\alpha_i\right)=-V_0(S\lambda)-V_0(V_0^i)\alpha_i\stackrel{(\ref{wag01}), (\ref{wag015})}{=}-V_0(S\lambda)-\lambda S(\lambda).$$
Finally we have
\begin{equation}
\label{wag05}
V_0^h (\lambda) =V_0(\lambda) V_0^i f_i \sqrt{g(V,V)}.
\end{equation}
Differentiating (\ref{wag04}) along the indicatrix curve
$$V_0\left(S \lambda\right)=[V_0,S](\lambda)+S\left(V_0 \lambda\right)\stackrel{(\ref{wag02})}{=} V_0^h (\lambda)+S\left(V_0 \lambda\right),$$
$$V_0 \left(V_0(\lambda) y^i f_i \sqrt{g(V,V)}\right)\stackrel{(\ref{difeq015})}=V_0 \left(V_0 \lambda \right) y^i f_i \sqrt{g(V,V)}+V_0(\lambda) V_0^i f_i \sqrt{g(V,V)}+V_0(\lambda) y^i \alpha_i\stackrel{(\ref{wag01})}{=}$$
$$V_0 \left(V_0 \lambda \right) y^i f_i \sqrt{g(V,V)}+V_0(\lambda) V_0^i f_i \sqrt{g(V,V)}$$
and, consequently,
$$V_0(\lambda)V_0^h (\lambda)+V_0(\lambda)S\left(V_0 \lambda\right)=V_0 \left(V_0 \lambda \right) V_0(\lambda) y^i f_i \sqrt{g(V,V)}+V_0(\lambda) V_0(\lambda)V_0^i f_i \sqrt{g(V,V)}\stackrel{(\ref{wag04}), (\ref{wag05})}{=}$$
$$V_0 \left(V_0 \lambda \right) S(\lambda)+V_0(\lambda) V_0^h(\lambda),$$
i.e.
\begin{equation}
\label{wag06}
V_0(\lambda)S\left(V_0 \lambda\right)=V_0 \left(V_0 \lambda \right) S(\lambda).
\end{equation}
In a similar way, differentiating (\ref{wag05}) along the indicatrix curve 
$$V_0 \left(V_0^h \lambda \right) =[V_0, V_0^h] \lambda +V_0^h \left(V_0 \lambda \right)\stackrel{(\ref{wag02})}{=}-\frac{1}{F}S_0 (\lambda)-\lambda V_0^h (\lambda) -S(\lambda) V_0 (\lambda)+V_0^h \left(V_0 \lambda \right),$$
$$V_0\left(V_0 (\lambda) V_0^i f_i \sqrt{g(V,V)}\right)\stackrel{(\ref{difeq015})}{=} V_0\left(V_0 \lambda \right) V_0^i f_i \sqrt{g(V,V)}+V_0(\lambda) V_0 \left(V_0^i\right)f_i \sqrt{g(V,V)}+V_0(\lambda)V_0^i \alpha_i\stackrel{(\ref{wag01}), (\ref{wag015})}{=}$$
$$ V_0\left(V_0 \lambda \right) V_0^i f_i \sqrt{g(V,V)}-\lambda V_0(\lambda) V_0^i f_i \sqrt{g(V,V)}-V_0(\lambda) \frac{y^i}{F^2}f_i \sqrt{g(V,V)}-V_0(\lambda) S(\lambda)$$
and, consequently,
$$-V_0(\lambda) \left(\frac{1}{F}S_0 (\lambda)+\lambda V_0^h (\lambda) +S(\lambda) V_0 (\lambda)-V_0^h \left(V_0 \lambda \right)\right)= $$
$$V_0\left(V_0 \lambda \right) V_0(\lambda) V_0^i f_i \sqrt{g(V,V)}-\lambda V_0(\lambda) V_0(\lambda) V_0^i f_i \sqrt{g(V,V)}-V_0(\lambda) V_0(\lambda) \frac{y^i}{F^2}f_i \sqrt{g(V,V)}-$$
$$V_0(\lambda) V_0(\lambda) S(\lambda)\stackrel{(\ref{wag04}), (\ref{wag05})}{=}V_0\left(V_0 \lambda \right) V_0^h(\lambda)-\lambda V_0(\lambda) V_0^h(\lambda) -\frac{1}{F}V_0(\lambda) S_0(\lambda)-V_0(\lambda) V_0(\lambda) S(\lambda),$$
i.e.
\begin{equation}
\label{wag07}
V_0(\lambda) V_0^h \left(V_0 \lambda \right)=V_0\left(V_0 \lambda \right) V_0^h(\lambda).
\end{equation}
Since $S$ and $V_0^h$ span the horizontal subspaces we can write, by  (\ref{wag06}) and (\ref{wag07}), that
\begin{equation}
\label{wagfinal}
V_0(\lambda) X_i^h\left(V_0 \lambda \right)=V_0\left(V_0 \lambda \right) X_i^h(\lambda) \ \ (i=1,2).
\end{equation}
Equations (\ref{wagfinal}) are called Wagner's equations \cite[Formula 18]{Wag1}. 
\vspace{0.2cm}
\begin{center}
\begin{tabular}{|l|l|l|}
\hline
 \multicolumn{3}{|c|}{}\\
 \multicolumn{3}{|c|}{Wagner's notations \cite{Wag1}}\\
\multicolumn{3}{|c|}{}\\
\hline
the evaluation of the main scalar & & the canonical horizontal \\
along the central affine arcwise&$\displaystyle{\frac{\partial A}{\partial \theta}}=(\lambda\circ c_p)'=V_0(\lambda)\circ c_p$&sections:\\
parametrization: $A=\lambda\circ c_p$&&$\displaystyle{\nabla_{\beta}=X_{\beta}^h}$, $\beta=1, 2$\\
\hline
\end{tabular}
\end{center}
\vspace{0.2cm}
Consider the indicatrix bundle $IM:=F^{-1}(1)$. Wagner's equations imply that 
\begin{equation}
\label{closedness}
V_0(\lambda) d \left(V_0 \lambda \right)=V_0\left(V_0 \lambda \right) d \lambda
\end{equation}
holds on the manifold $IM$ because $\displaystyle{V_0(\lambda) V_0 \left(V_0 \lambda \right)=V_0\left(V_0 \lambda \right) V_0 (\lambda)}$ is automathic; note that
$$V_0(F)=X_i^h(F)=0 \ \ (i=1, 2),$$
i.e. $V_0$, $X_1^h$ and $X_2^h$ form a local frame of the indicatrix bundle. Suppose that $F(v)=1$ and $V_0 (v) \lambda \neq 0$. Equation (\ref{closedness}) implies that $\displaystyle{d \left(V_0 \lambda \right)}$ is the proportional of $d\lambda$ around $v$ and, consequently,
$$d(V_0\lambda) \wedge d\lambda=0 \ \ \Leftrightarrow \ \ d \left((V_0\lambda) d\lambda\right)=0.$$
This means that there is a (local) solution $\mu$ such that 
\begin{equation}
\label{wag08}
(V_0\lambda) d\lambda=d \mu.
\end{equation}
Taking a coordinate system $\varphi=(z^1, z^2, \lambda)$ of the indicatrix bundle around $v$, formula (\ref{wag08}) says that $\displaystyle{\partial \mu/\partial z^1=\partial \mu/\partial z^2=0}$. This means that $\mu$ depends only on $\lambda$. If the function $f$ is defined by $\displaystyle{f(\lambda):=\mu'(\lambda )}$, where $\mu$ is the local solution of (\ref{wag08}), then $\displaystyle{V_0 (\lambda)=f(\lambda)}$ as Wagner's theorem states.

\subsubsection{Wagner's theorem} \cite{Wag1} \emph{A necessary and sufficient condition that $\displaystyle{F_2\ \left(\frac{\partial A}{\partial \theta}\neq 0\right)}$ be a generalized Berwald space is that $\displaystyle{\frac{\partial A}{\partial \theta}}$ be a function of $A$}.

\begin{Exc}
Prove the converse of Wagner's theorem.
\end{Exc} 

\section{The averaging method} Let $\nabla$ be a linear connection on the base manifold $M$ of dimension $2$ and suppose that the parallel transports preserve the Finslerian length of tangent vectors (compatibility condition). By the fundamental result of the theory \cite{V5} such a linear connection must be metrical with respect to the averaged Riemannian metric given by integration of the Riemann-Finsler metric on the indicatrix hypersurfaces; see Theorem \ref{heritage}. Therefore it is uniquely determined by its torsion tensor of the form
\begin{equation}
\label{torsionform}
T(X,Y)=\rho(X)Y-\rho(Y)X;
\end{equation}
see Formula \ref{tor}. The idea of the comparison of $\nabla$ with the L\'{e}vi-Civita connection $\nabla^*$ associated with the averaged Riemannian metric (\ref{averagemetric1}) was used to solve the problem of the intrinsic characterization of the semi-symmetric compatible linear connections for both low and higher dimensional spaces. The solution is the expression of the $1$ - form $\rho$ in terms of metrics and differential forms given by averaging. For the details see  \cite{V6}, \cite{V9}, \cite{V10} and \cite{V11}. 

\subsection{The divergence representation of the Gauss curvature} Let a point $p\in M$ be given and consider the orthogonal group with respect to the averaged Riemannian metric. It is clear that the subgroup $G\subset O(2)$ of the orthogonal transformations leaving the Finslerian indicatrix invariant is finite unless the Finsler surface reduces to a Riemannian one; see also \cite{V12}. If $\nabla$ is a linear connection on the base manifold such that the parallel transports preserve the Finslerian length of tangent vectors (compatibility condition) then $\textrm{Hol}_p \nabla \subset G$ is also finite for any $p\in M$ and the curvature tensor of $\nabla$ is zero. In what follows we are going to compute the relation between the curvatures of $\nabla$ and $\nabla^*$. Taking vector fields with pairwise vanishing Lie brackets on the neighbourhood $U$ of the base manifold, the Christoffel process implies that
\begin{equation}
\label{Cproc}
\gamma(\nabla^*_X Y,Z)=\gamma(\nabla_X Y, Z)+\frac{1}{2}\left(\gamma(X, T(Y,Z))+\gamma(Y, T(X,Z))-\gamma(Z, T(X,Y))\right),
\end{equation}
where $\nabla^*$ denotes the L\'{e}vi-Civita connection. If the torsion is of the form (\ref{torsionform}), then we have that
\begin{eqnarray}
\nabla^*_X Y =  \nabla_X Y +  \rho (Y) X  -  \gamma ( X, Y ) \rho^{\sharp} \ \ \Rightarrow \ \ \nabla_X Y =  \nabla^*_X Y -  \rho (Y) X  +  \gamma ( X, Y ) \rho^{\sharp}, \label{02}
\end{eqnarray}
where $\rho^{\sharp}$ is the dual vector field of $\rho$ defined by $\displaystyle{\gamma(\rho^{\sharp}, X)=\rho(X)}$. 
Consider the curvature tensor
\begin{eqnarray}
R(X,Y)Z= \nabla_{X} \nabla_{Y} Z - \nabla_{Y} \nabla_{X} Z \label{03}
\end{eqnarray}
of $\nabla$. Substituting (\ref{02}) into (\ref{03})
\begin{eqnarray*}
R(X,Y)Z= \nabla_{X} \left(\nabla^*_Y Z -  \rho (Z) Y  + \rho^{\sharp} \gamma ( Z, Y ) \right)-\nabla_{Y} \left(\nabla^*_X Z -  \rho (Z) X  + \rho^{\sharp} \gamma ( Z, X ) \right).
\end{eqnarray*}
Some further direct computations show that 
$$R(X,Y)Z= R^*(X,Y)Z+$$
$$\bigg(\gamma(X,Z) \|\rho^{\sharp}\|^2-\rho(X)\rho(Z)-\left(\nabla_X^* \rho\right)(Z)\bigg)Y+\gamma(Y,Z)\nabla^*_X \rho^{\sharp}+\gamma(Y,Z)\rho(X)\rho^{\sharp}+$$
$$\bigg(\left(\nabla_Y^* \rho \right)(Z)+\rho(Y)\rho(Z)-\gamma(Y,Z) \|\rho^{\sharp}\|^2\bigg) X-\gamma(X,Z)\nabla^*_Y \rho^{\sharp}-\gamma(X,Z)\rho(Y)\rho^{\sharp}.$$
Since the holonomy group of $\nabla$ must be finite in case of a non-Riemannian generalized Berwald surface we have that $R(X,Y)Z=0$. Taking an orthonormal frame $\gamma(X,Y)=0$, $\gamma(X,X)=\gamma(Y,Y)=1$ at the point of $p\in M$ it follows that 
$$0=\gamma (R^*(X,Y)Y,X)+\rho^2(X)+\rho^2(Y)-\|\rho^{\sharp}\|^2+\gamma \left(\nabla^*_X \rho^{\sharp},X \right)+\left(\nabla^*_Y \rho \right)(Y),$$
where $\displaystyle{\rho^2(X)+\rho^2(Y)-\|\rho^{\sharp}\|^2=0}$ and
$$\left(\nabla^*_Y \rho \right)(Y)=Y\rho(Y)-\rho(\nabla^*_Y Y)=Y\gamma(\rho^{\sharp}, Y)-\rho(\nabla^*_Y Y)=\gamma(\nabla^*_Y \rho^{\sharp}, Y).$$
Therefore
\begin{equation}
\label{cur}
0=\kappa^*(p)+\textrm{div}^* \rho^{\sharp}(p) \ \ \Rightarrow \ \ \kappa^*=-\textrm{div}^* \rho^{\sharp},
\end{equation}
where $\kappa^*$ is the Gauss curvature of the manifold with respect to the averaged Riemannian metric and $\displaystyle{\textrm{div}^* \rho^{\sharp}:=\gamma \left(\nabla^*_X \rho^{\sharp},X \right)+\gamma \left(\nabla^*_Y \rho^{\sharp},Y \right)}$
is the divergence operator. Equation (\ref{cur}) is called \emph{the divergence representation of the Gauss curvature}.

\begin{Cor} A Riemannian surface admits a metric linear connection of zero curvature if and only if its Gauss curvature can be represented as a divergence of a vector field.
\end{Cor}

\begin{Cor} If $M$ is a compact generalized Berwald surface without boundary then it must have zero Euler characteristic.
\end{Cor}
\begin{Pf}
Taking the integral of the divergence representation (\ref{cur}) we have the zero Euler characteristic due to the Gauss-Bonnet theorem and the divergence theorem.
\end{Pf}

\begin{Cor}
\label{forbidden} A two-dimensional Euclidean sphere could not carry Finslerian structures admitting compatible linear connections. 
\end{Cor}

\begin{Rem}
Using the classification of orientable compact surfaces without boundary we can also state that they could not carry Finslerian structures admitting compatible linear connections except the case of genus $1$. 
\end{Rem}

\subsection{Exact and closed Wagner manifolds} It is well-known that any Riemannian surface is locally conformally flat by the (local) solution of the second order elliptic partial differential equation $\displaystyle{\Delta^* f=\kappa^*}$. Its Finslerian analogue is that \emph{any non-Riemannian Finsler surface is locally conformal to a locally Minkowski manifold of dimension $2$}; a locally Minkowski manifold is a Berwald manifold (torsion-free case, i.e. the compatible linear connection is $\nabla^*$) such that $R^*=0$. The solution of the so-called Matsumoto's problem \cite{V6}, see also \cite{V7}, proves that the statement is false in the non-Riemannian Finsler geometry. 
\begin{itemize}
\item[Step 1] By Hashiguchi and Ichyjio's classical theorem \cite{HY}, see also \cite{V1} and \cite{V2}, a Finsler manifold is a conformally Berwald manifold if and only if there exists a semi-symmetric compatible linear connection  with an exact $1$-form $\rho$ in the torsion (\ref{torsionform}). Especially, it is the exterior derivative of the logarithmic scale function $\alpha$ between the (conformally related, see \cite{H2}) Finslerian fundamental functions $\displaystyle{\tilde{F}=e^{\alpha\circ \pi}F}$ up to a minus sign. 
\end{itemize}
\begin{Def} Generalized Berwald manifolds admitting compatible semi-symmetric linear connections with an exact $1$-form $\rho$ in the torsion \eqref{torsionform} are called exact Wagner manifolds. Generalized Berwald manifolds admitting compatible semi-symmetric linear connections with a closed $1$-form $\rho$ in the torsion \eqref{torsionform} are called closed Wagner manifolds.
\end{Def}

\begin{itemize}
\item[Step 2] The generalization of Hashiguchi and Ichyjio's classical theorem for closed Wagner manifolds is the statement that \emph{a Finsler manifold is a locally conformally Berwald manifold if and only if it is a closed Wagner manifold}. It is clear from the global version of the theorem that any point of a closed Wagner manifold has a neighbourhood over which it is conformally equivalent to a Berwald manifold, i.e. any closed Wagner manifold is a locally conformally Berwald manifold. What about the converse? Suppose that we have a locally conformally Berwald manifold. The exterior derivatives of the local scale functions constitute  a globally well-defined closed $1$-form for the torsion (\ref{torsionform}) of a compatible linear connection if and only if they coincide on the intersections of overlapping neighbourhoods. Since the conformal equivalence is transitive it follows that overlapping neighbourhoods carry conformally equivalent Berwald metrics. The problem posed by M. Matsumoto \cite{M2} in 2001 is that \emph{are there non-homothetic and non-Riemannian conformally equivalent Berwald spaces}? It has been completely solved by \cite{V6} in 2005, see also \cite{V7}.
\end{itemize}

\begin{Thm} \emph{\cite{V6}, see also \cite{V7}} The scale function between conformally equivalent Berwald manifolds must be constant unless they are Riemannian. 
\end{Thm}

\begin{Cor} \emph{\cite{V6}, see also \cite{V7}} A Finsler manifold is a locally conformally Berwald manifold if and only if it is a closed Wagner manifold.
\end{Cor}

Using Corollary \ref{forbidden} we have the following result.

\begin{Cor} A two-dimensional Euclidean sphere could not carry non-Riemannian locally conformally Berwald Finslerian structures. Especially, it can not be a locally conformally flat non-Riemmanian Finsler manifold.
\end{Cor}

\subsection{The case of the Euclidean plane} Consider the Euclidean plane $\mathbb{R}^2$ equipped with the canonical inner product $\delta_{ij}$ and let $\displaystyle{\rho=\rho_1du^1+\rho_2du^2}$ be a $1$-form. If $X$ is a parallel vector field with respect to $\nabla$ along a curve $\displaystyle{c\colon [0, 1] \to \mathbb{R}^2}$, then, by formula (\ref{02}),
$$\left(X^1\right)'=\left(c^1\right)'\rho_c (X)-\left(\left(c^1\right)'X^1+\left(c^2\right)'X^2\right)\rho^1\circ c,$$
$$\left(X^2\right)'=\left(c^2\right)'\rho_c (X)-\left(\left(c^1\right)'X^1+\left(c^2\right)'X^2\right)\rho^2\circ c,$$
where $\rho^i=\delta^{ik}\rho_k=\rho_i$ ($i=1, 2$). Therefore
$$\left(X^1\right)'=X^2\left(\left(c^1\right)'\rho_2 \circ c-\left(c^2\right)'\rho_1\circ c\right) \ \ \textrm{and}\ \ \left(X^2\right)'=X^1\left(\left(c^2\right)'\rho_1 \circ c-\left(c^1\right)'\rho_2\circ c \right).$$
If the divergence of $\rho^{\sharp}$ vanishes (the curvature of the Euclidean plane is identically zero), then the curl of its rotated vector field
$$\rho_2 \frac{\partial}{ \partial u^1}-\rho_1 \frac{\partial}{ \partial u^2}$$
is zero, i.e. we have a global solution (potential) of equations $\displaystyle{\rho_2=\frac{\partial f}{ \partial u^1}}$ and $\displaystyle{\rho_1=-\frac{\partial f}{ \partial u^2}}$. Therefore the differential equations of the parallel vector fields are $\displaystyle{\left(X^1\right)'=\varphi' X^2}$ and $\displaystyle{\left(X^2\right)'=-\varphi' X^1}$, where $\displaystyle{\varphi=f\circ c}$. Since $\nabla$ is metrical 
$$X(t)=r_0\left(\cos \theta(t) \frac{\partial}{ \partial u^1}\circ c(t)+\sin \theta(t) \frac{\partial}{ \partial u^2}\circ c(t)\right)$$
with constant Euclidean norm $r_0$, where $\theta'(t)=-\varphi'(t)$ because of the parallelism. The general form of a parallel vector field with respect to $\nabla$ along the curve $\displaystyle{c\colon [0, 1] \to \mathbb{R}^2}$ is
$$X(t)=r_0\left(\cos \left(\varphi(t)+\varphi_0\right) \frac{\partial}{ \partial u^1}\circ c(t)-\sin \left(\varphi(t)+\varphi_0\right) \frac{\partial}{ \partial u^2}\circ c(t)\right).$$
It is clear that if $c(0)=c(1)$, then $X(0)=X(1)$, i.e. the holonomy group of $\nabla$ contains only the identity. Taking an arbitrary convex curve around the origin we can extend it by parallel transports with respect to $\nabla$ to the entire plane $\mathbb{R}^2$. Such a collection of indicatrices constitutes a Finslerian metric function $F$ with $\nabla$ as the compatible linear connection.

\subsubsection{An example} If $\displaystyle{\rho=u^2du^1-u^1du^2}$ then $\displaystyle{f(u^1, u^2)=-\frac{1}{2}\left(\left(u^1\right)^2+\left(u^2\right)^2\right)}$ and the parallel vector fields are of the form
$$X(t)=X^1(t)\frac{\partial}{ \partial u^1}\circ c(t)+X^2(t)\frac{\partial}{ \partial u^2}\circ c(t),$$
where
$$X^1(t)=r_0 \cos \left(\frac{1}{2}\left(\left(c^1\right)^2(t)+\left(c^2\right)^2(t)\right)+\varphi_0\right),$$
$$X^2(t)=-r_0\sin \left(\frac{1}{2}\left(\left(c^1\right)^2(t)+\left(c^2\right)^2(t)\right)+\varphi_0\right).$$
Let the trifocal ellipse defined by 
\begin{equation}
\label{3ell}
\sqrt{(u^1+1)^2+(u^2)^2}+\sqrt{(u^1)^2+(u^2)^2}+\sqrt{(u^1-1)^2+(u^2)^2}=4
\end{equation}
be choosen as the indicatrix at the origin. The focal set contains the elements
$$-X_0:=(-1,0), \ {\bf 0}, \ X_0:=(1,0).$$
The parallel translates of the trifocal ellipse (\ref{3ell}) are given by the equations
\begin{equation}
\label{3elltrans}
\sqrt{(u^1+X^1(t))^2+(u^2+X^2(t))^2}+\sqrt{(u^1)^2+(u^2)^2}+\sqrt{(u^1-X^1(t))^2+(u^2-X^2(t))^2}=4,
\end{equation}
where $X$ is a parallel vector field along a curve $c$ satisfying the initial conditions $X^1(0)=1$ and $X^2(0)=0$. The focal set at the parameter $t$ is $\displaystyle{-X(t), \ {\bf 0}, \ X(t)}$. Figure 1 shows the parallel translates of the indicatrix along the radial direction $c(t)=(t,t)$.
\begin{figure}[h]
\includegraphics[scale=0.25]{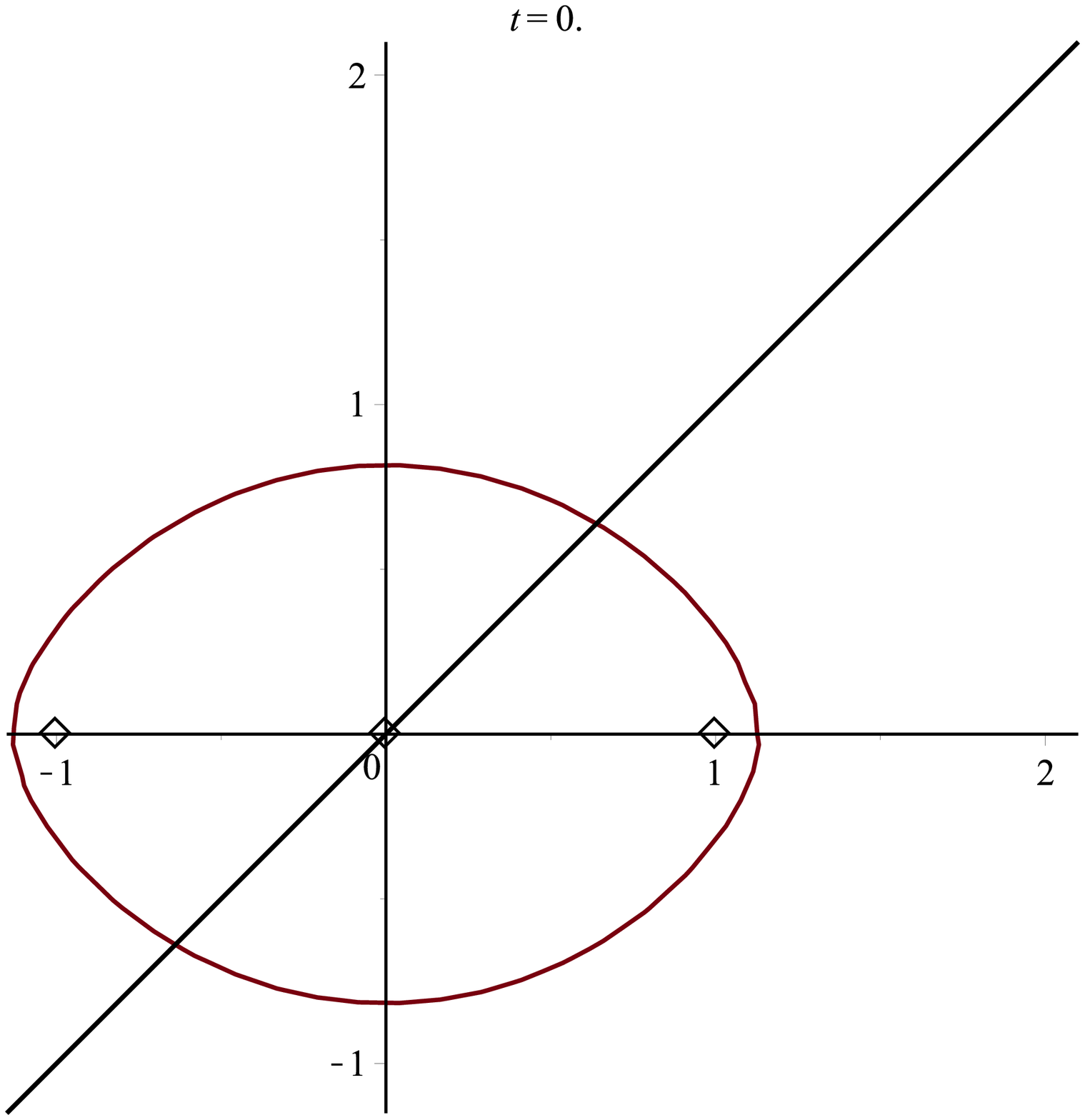}
\includegraphics[scale=0.25]{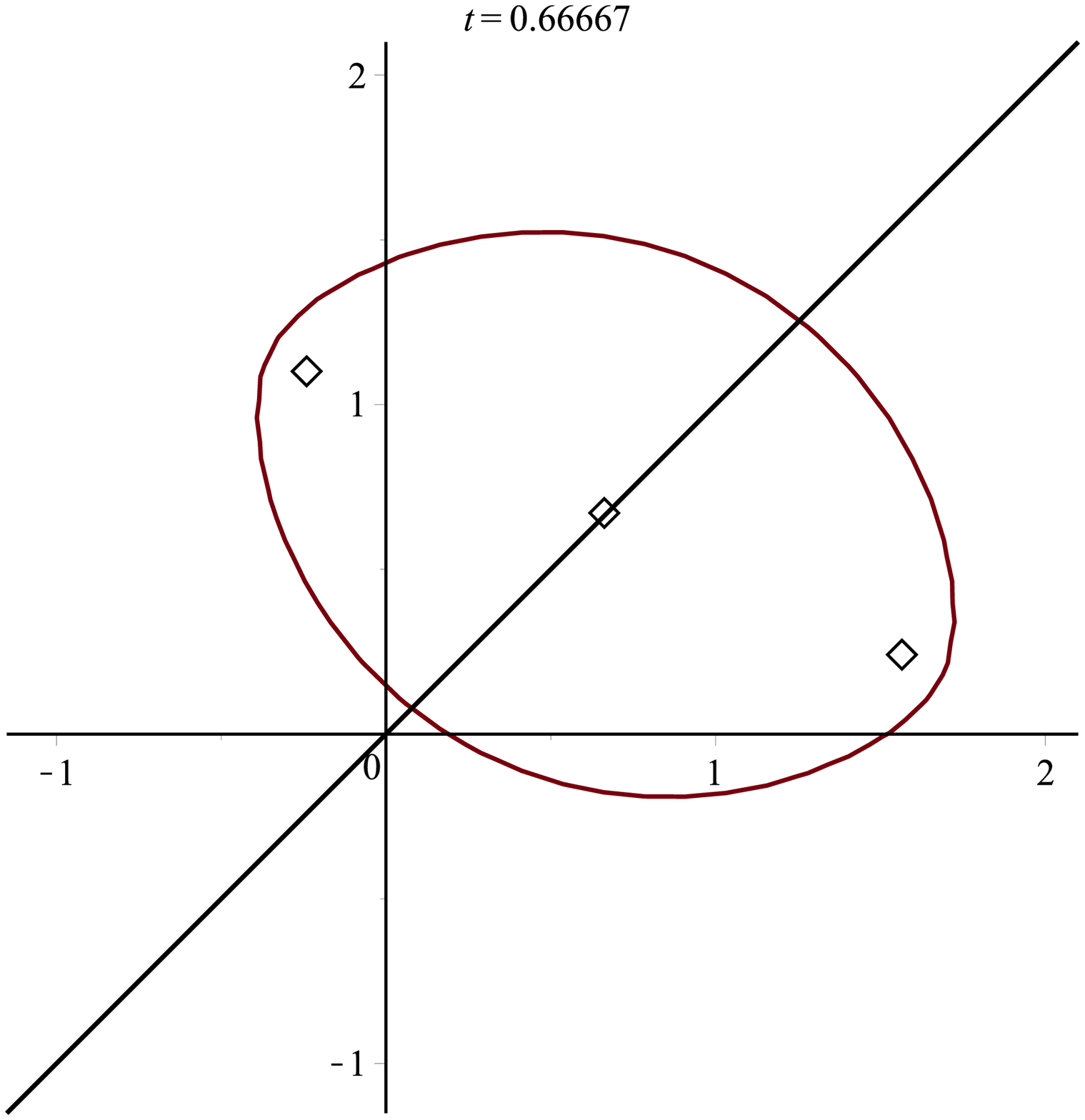}
\includegraphics[scale=0.25]{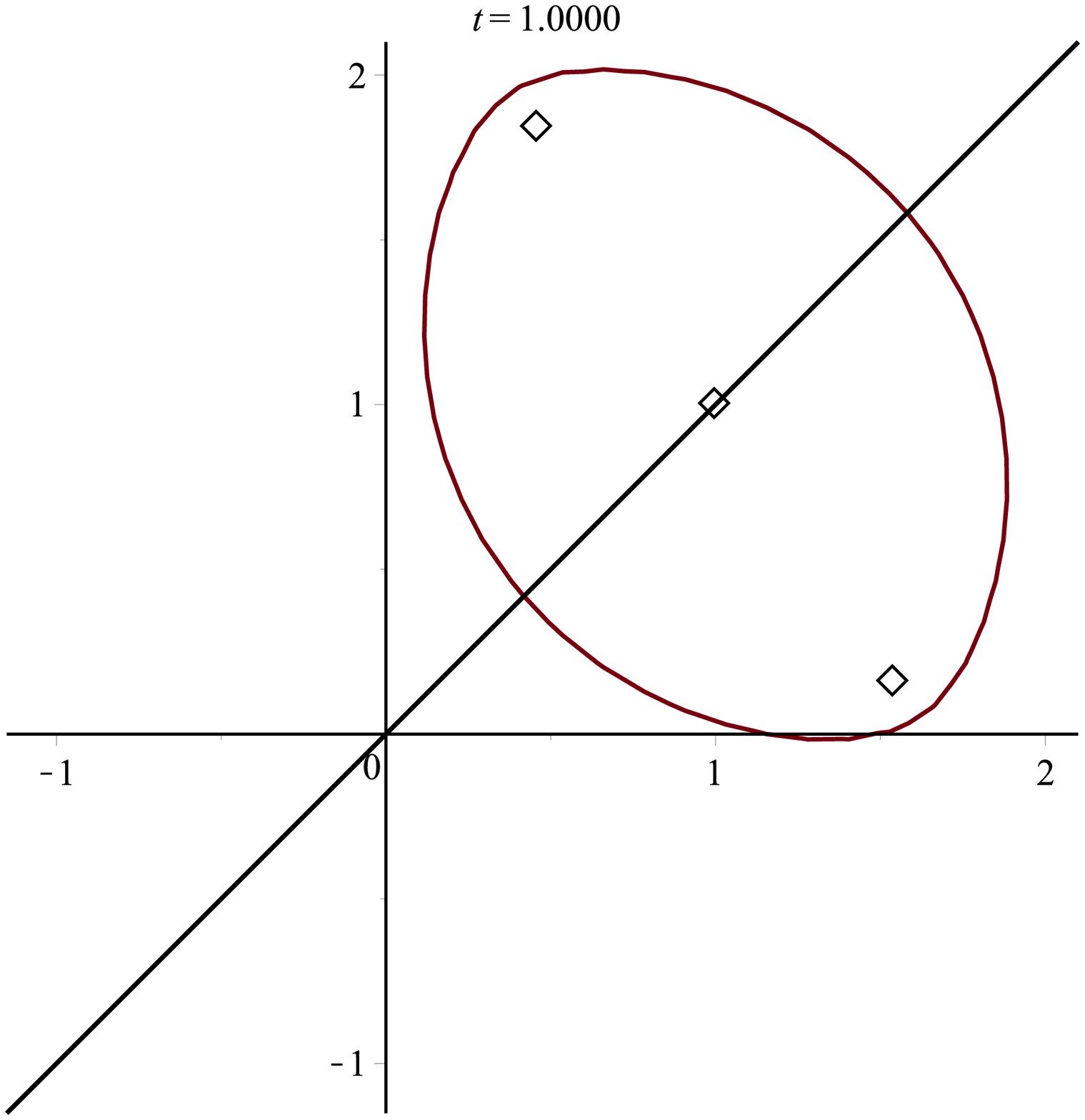}
\caption{Parallel translation along the radial direction (from left to right)}
\end{figure}
In case of the radial direction the focal set of the translated trifocal ellipse at the parameter $t$ is $-X(t)$, ${\bf 0}$,
$$X(t)=\cos \left(t^2\right) \frac{\partial}{ \partial u^1}_{(t,t)}-\sin \left(t^2\right) \frac{\partial}{ \partial u^2}_{(t,t)}.$$

Figure 2 shows the parallel translates of the indicatrix along the circle $c(t)=(\cos(t), \sin(t)+1)$. The focal set of the translated trifocal ellipse at the parameter $t$ is $-X(t)$, ${\bf 0}$, 
$$X(t)=\cos \left(1+\sin(t)\right) \frac{\partial}{ \partial u^1}_{(\cos(t),1+\sin(t))}-\sin \left(1+\sin(t)\right) \frac{\partial}{ \partial u^2}_{(\cos(t),1+\sin(t))}.$$
\begin{figure}[h]
\includegraphics[scale=0.25]{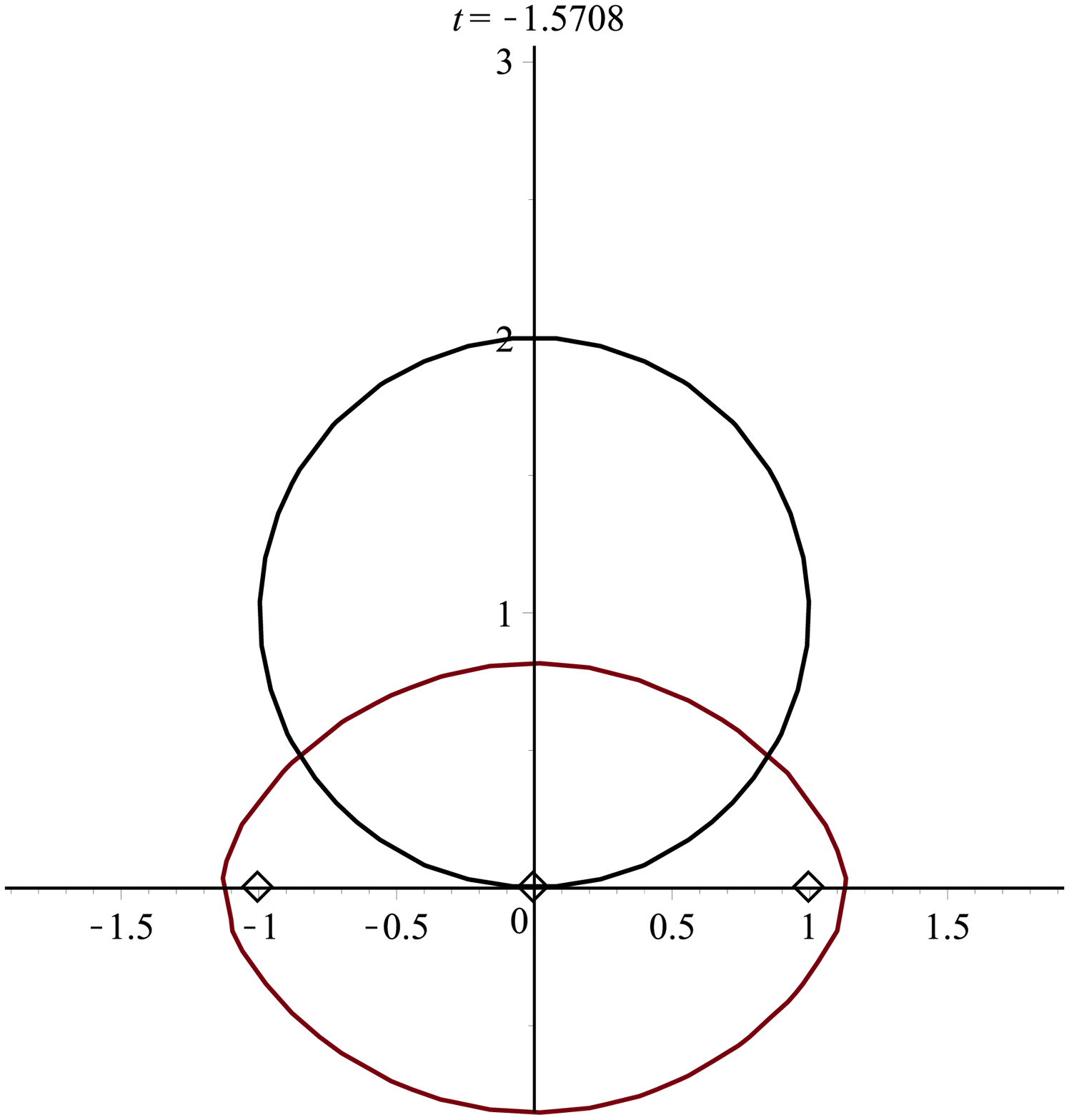}
\includegraphics[scale=0.25]{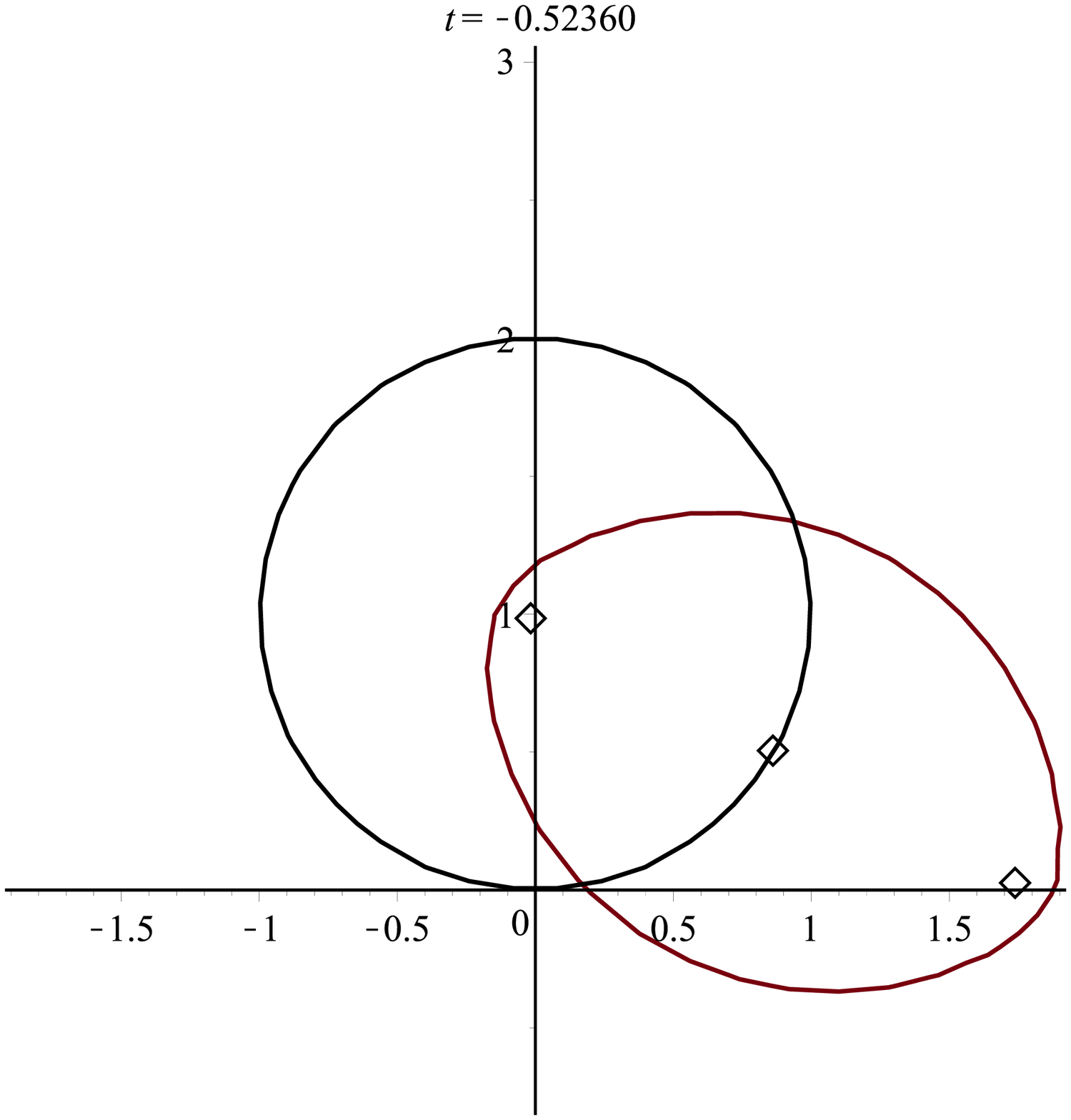}
\includegraphics[scale=0.25]{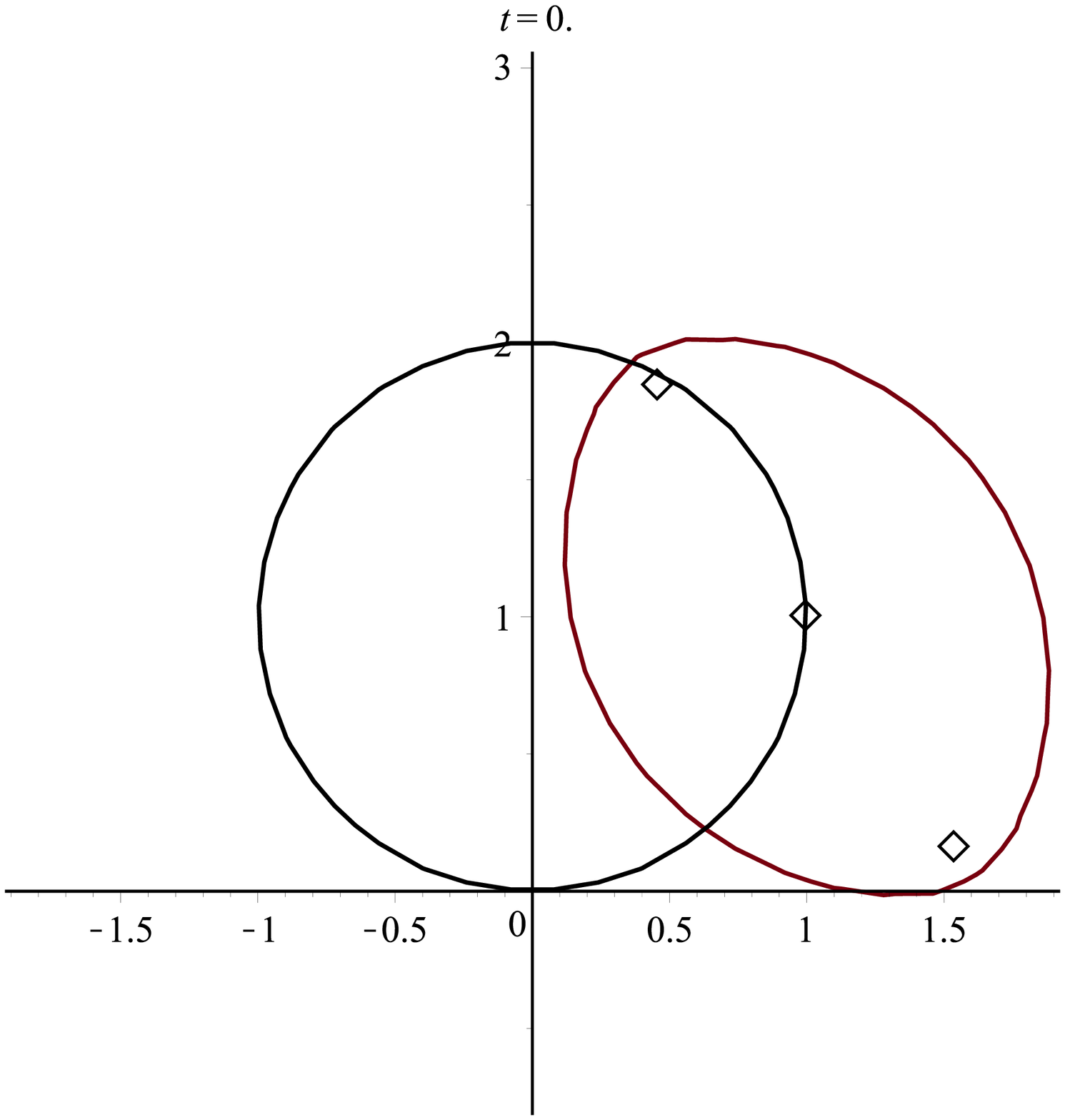}
\end{figure}
\begin{figure}[h]
\includegraphics[scale=0.25]{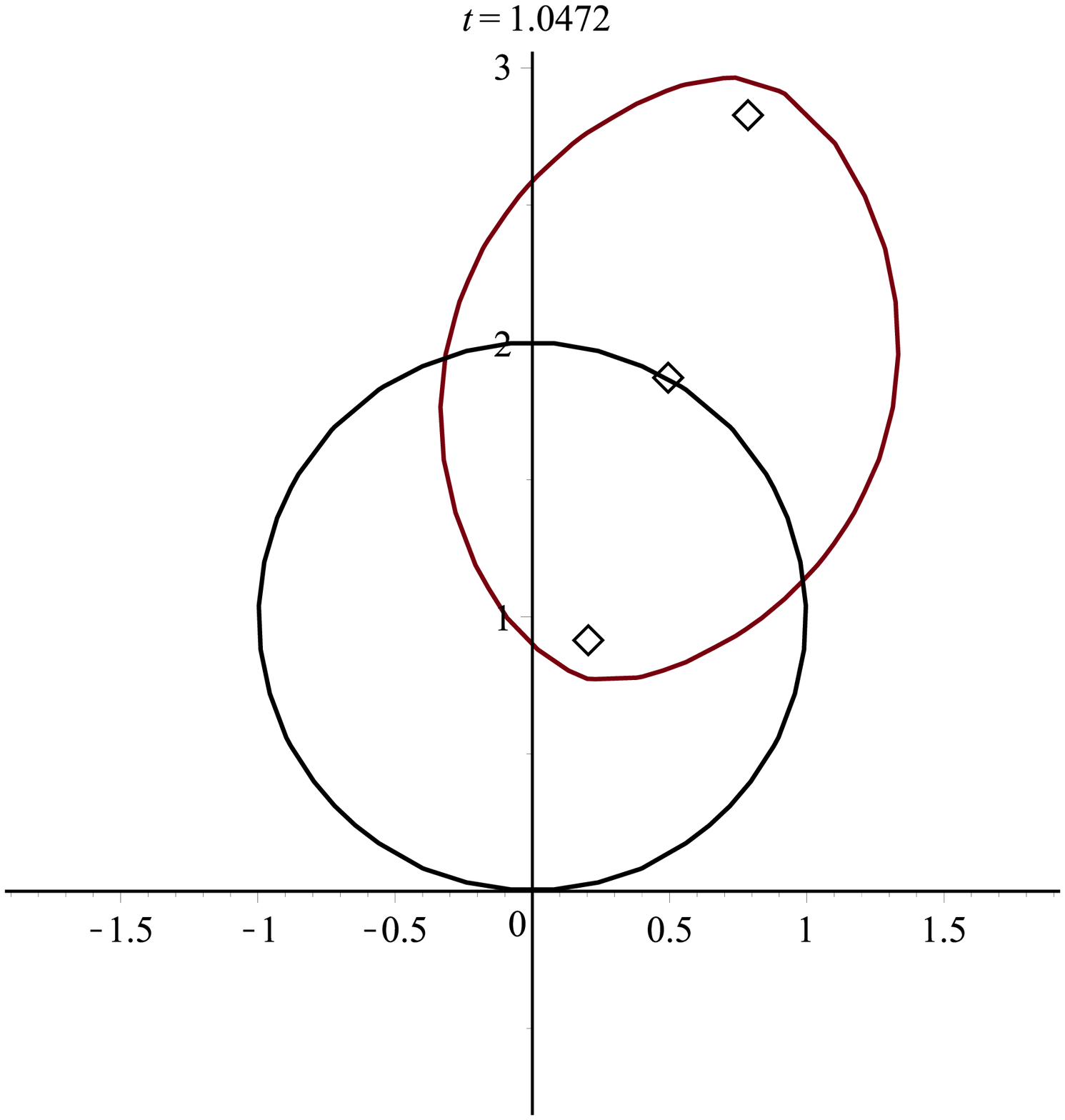}
\includegraphics[scale=0.25]{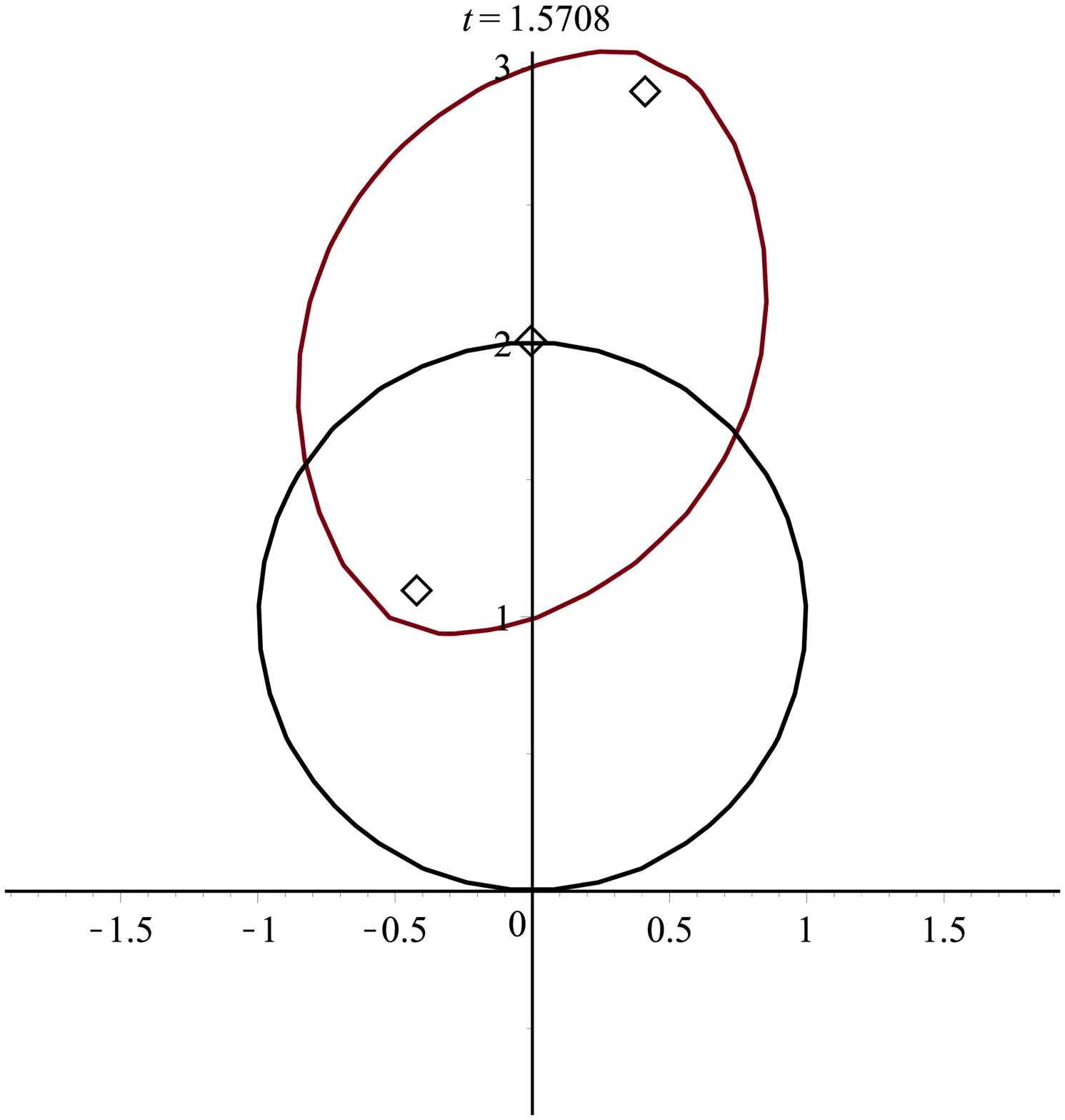}
\includegraphics[scale=0.25]{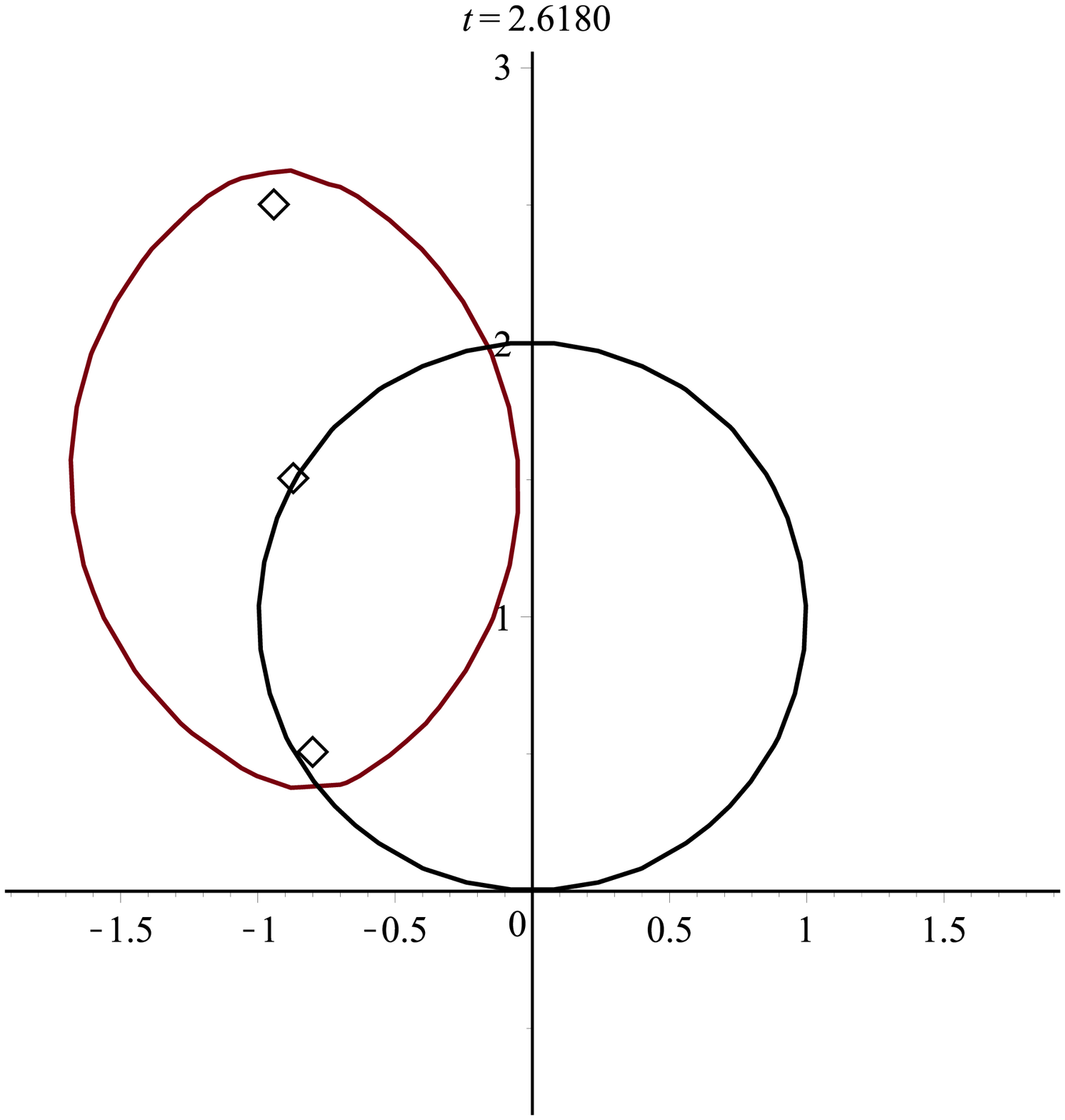}
\end{figure}
\begin{figure}[h]
\includegraphics[scale=0.25]{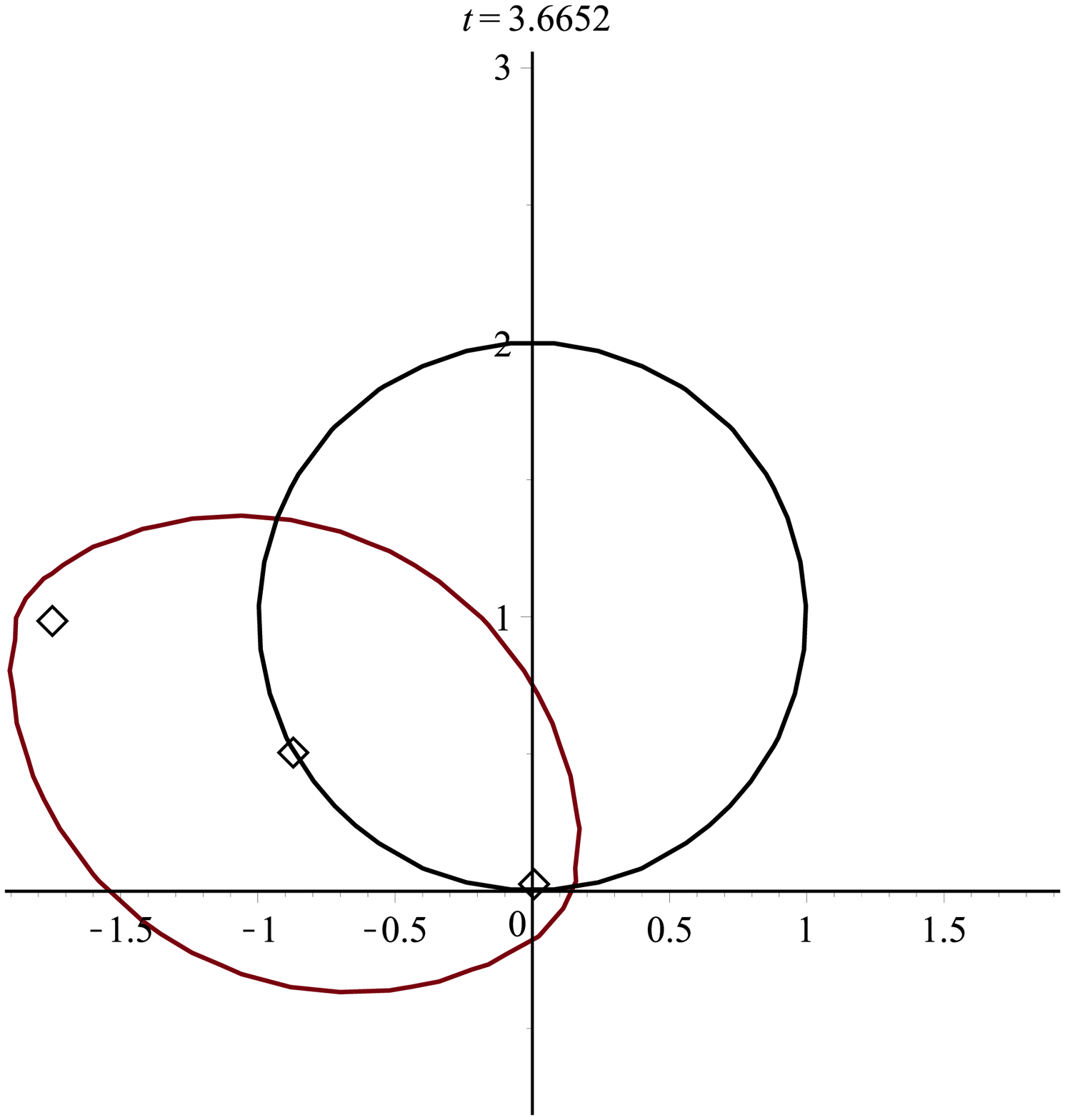}
\includegraphics[scale=0.25]{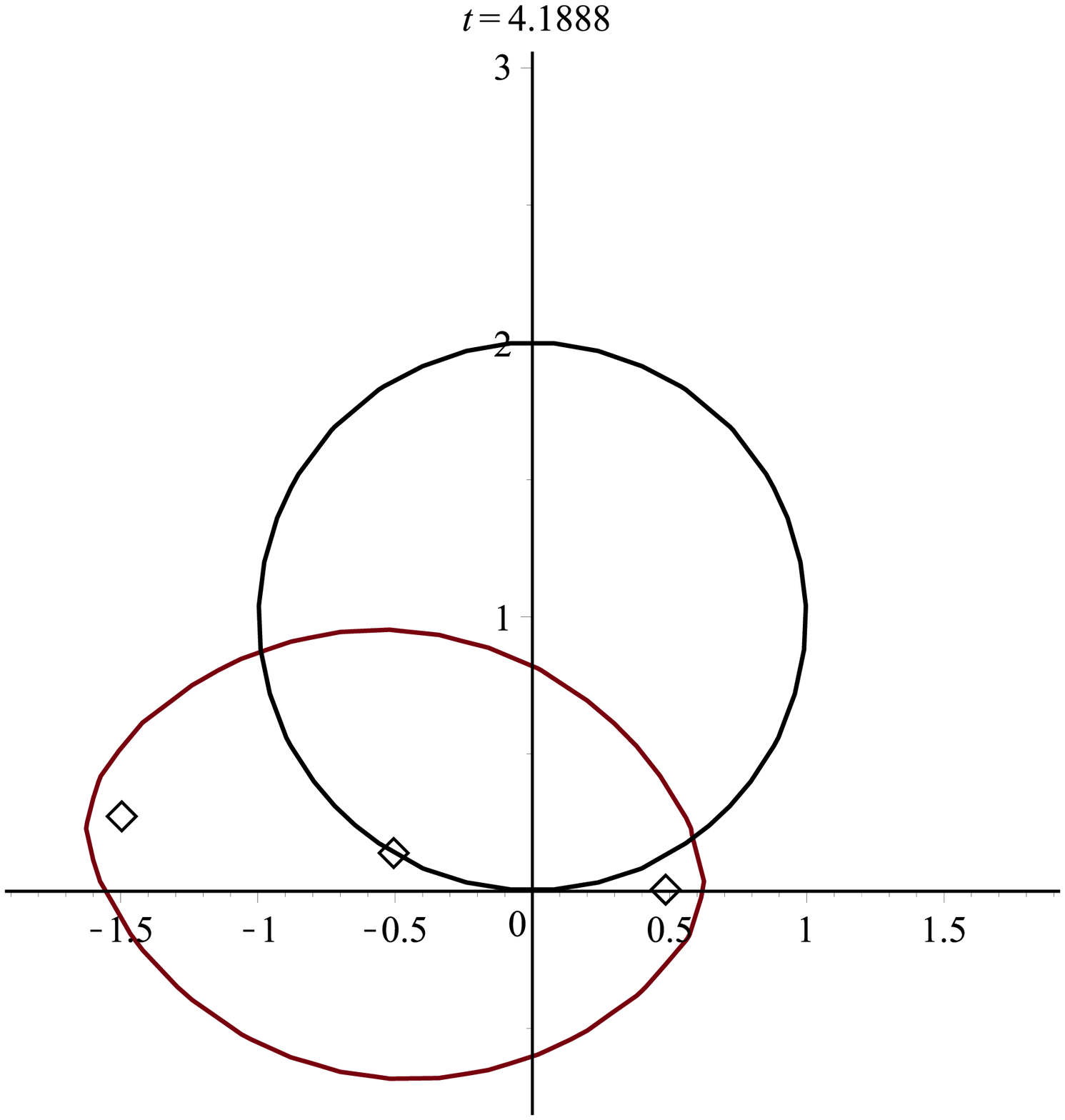}
\includegraphics[scale=0.25]{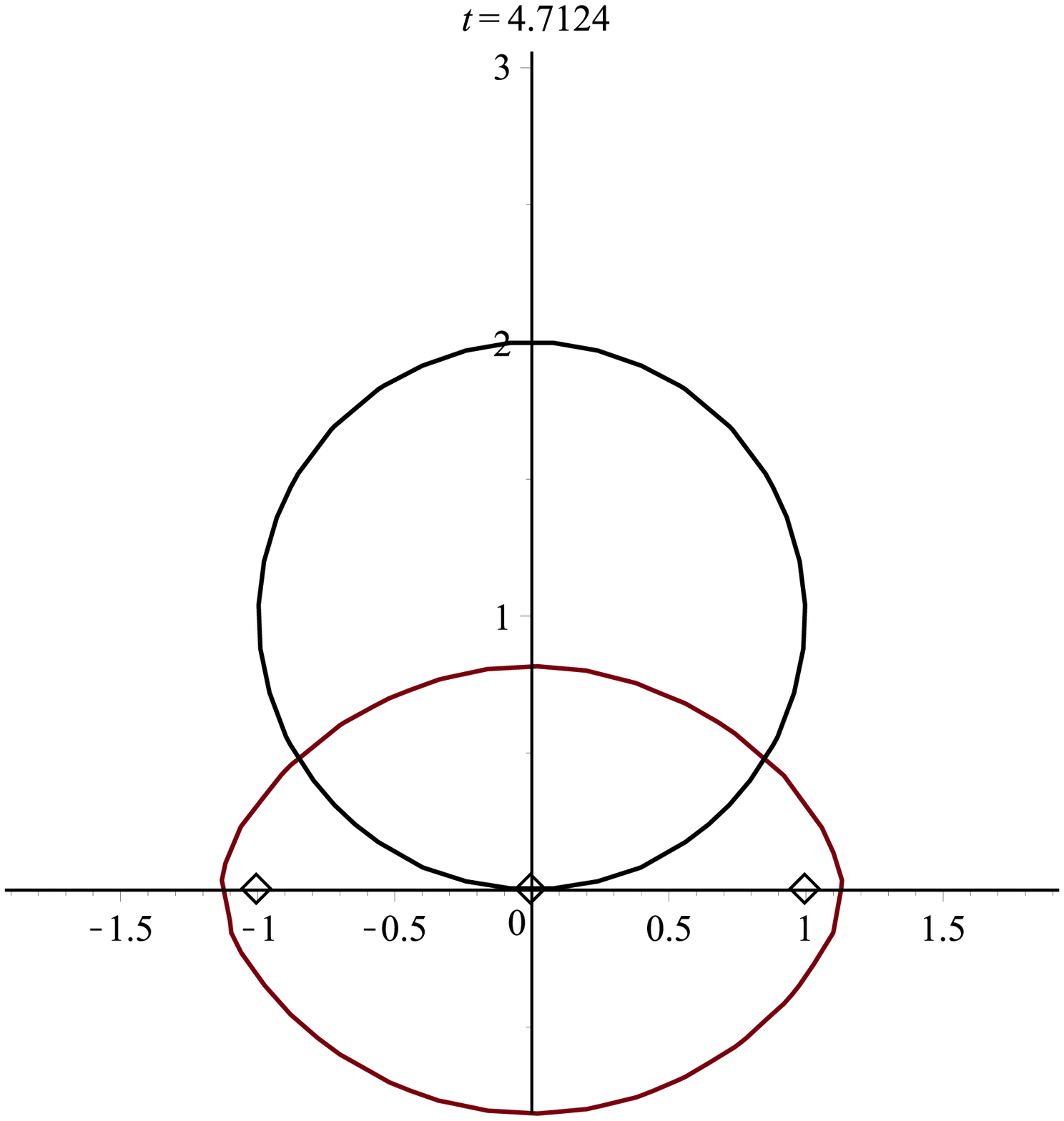}
\caption{Parallel translation along the circle (from left to right)}
\end{figure}

The induced generalized Berwald plane is not conformally flat. To present examples for conformally flat generalized Berwald manifold it is sufficient and necessary to choose a closed (and, consequently, exact) $1$-form $\rho$.

\end{document}